%% file: PflPosTanAITO.tex
\documentclass[reqno]{amsart}
\usepackage{amssymb,amsmath,amsfonts,latexsym}
\usepackage{pb-diagram}
\numberwithin{equation}{section}

\newcommand{\C}{{\mathbb{C}}}
\newcommand{\complex}{{\mathbb{C}}}
\newcommand{\R}{{\mathbb{R}}}

\newcommand{\Z}{{\mathbb{Z}}}

\newcommand{\N}{{\mathbb{N}}}

\newcommand{\K}{{\mathbb{K}}}
\newcommand{\W}{{\mathbb{W}}}
\newcommand{\A}{{\mathbb{A}}}

\newcommand{\frakg}{\mathfrak{g}}
\newcommand{\g}{\mathfrak{g}}

\newcommand{\frakh}{\mathfrak{h}}
\newcommand{\h}{\mathfrak{h}}

\newcommand{\mat}{\mathfrak{M}}

\newcommand{\calA}{\mathcal{A}}

\newcommand{\calC}{\mathcal{C}}

\newcommand{\calO}{\mathcal{O}}

\newcommand{\calS}{\mathcal{S}}

\newcommand{\Tr}{\operatorname{Tr}}
\newcommand{\tr}{\operatorname{tr}}

\newcommand{\Hom}{\operatorname{Hom}}
\newcommand{\End}{\operatorname{End}}
\newcommand{\im}{\operatorname{im}}

\newcommand{\gl}{\mathfrak{gl}}
\newcommand{\Sp}{\operatorname{Sp}}

\renewcommand{\det}{\operatorname{det}}
\newcommand{\Ch}{\operatorname{Ch}}
\newcommand{\ev}{\operatorname{ev}}
\newcommand{\supp}{\operatorname{supp}}
\newcommand{\Sym}{\operatorname{S}}
\newcommand{\ASym}{\operatorname{AS}}
\newcommand{\JSym}{\operatorname{JS}}
\newcommand{\PDO}{\operatorname{\Psi}}
\newcommand{\Op}{\operatorname{Op}}
\newcommand{\pr}{\operatorname{pr}}

\newcommand{\ad}{\operatorname{ad}}
\newcommand{\Exp}{\operatorname{Exp}}

\theoremstyle{plain}
        \newtheorem{theorem}{Theorem}[section]
        \newtheorem{lemma}[theorem]{Lemma}
        \newtheorem{proposition}[theorem]{Proposition}
        \newtheorem{corollary}[theorem]{Corollary}
        
\theoremstyle{definition}
        \newtheorem{definition}[theorem]{Definition}
        \newtheorem{remark}[theorem]{Remark}

\title{An algebraic index theorem for orbifolds}
\date{\today}
\author{M.J.~Pflaum, H.B.~Posthuma \textrm{and} X.~Tang}
\begin{document}
\begin{abstract}
 Using the concept of a twisted trace density on a cyclic groupoid, a 
 trace is constructed on a formal deformation quantization of a symplectic
 orbifold. An algebraic index theorem for orbifolds follows as a 
 consequence of a local Riemann--Roch theorem for such densities. 
 In the case of a reduced orbifold, this proves a conjecture by Fedosov, 
 Schulze, and Tarkhanov. Finally, it is shown how the Kawasaki index 
 theorem for elliptic operators on orbifolds follows from this algebraic 
 index theorem. 
\end{abstract}
\address{\newline
   Markus J. Pflaum, {\tt pflaum@math.uni-frankfurt.de}\newline
   \indent {\rm Fachbereich Mathematik, Goethe-Universit\"at Frankfurt/Main, 
           Germany } \newline
   Hessel Posthuma, {\tt posthuma@maths.ox.ac.uk}\newline
   \indent {\rm Mathematical Institute, University of Oxford, UK } \newline
   Xiang Tang, {\tt xtang@math.udavis.edu}   \newline   
   \indent {\rm  Department of Mathematics, University of California, Davis, 
           USA } 
}
\maketitle
\tableofcontents
\input{intro}
\input{preli}
\input{twtr}
\input{twhoch}
\input{cotwtrd}
\input{locrr}
\input{algind}

\appendix

\input{appa}
\input{appb}

\bibliographystyle{alpha}

\end{document}

%% file: intro.tex
\section*{Introduction}
Index theory originated with the seminal paper \cite{AtiSin} of Atiyah and
Singer almost 40 year ago. They proved that the index of an elliptic operator
on a closed Riemannian manifold $M$ depends only on the class of the principal 
symbol in the $K$-theory of $T^*M$. Ever since, many new
proofs and generalizations of this theorem have appeared. To mention, the
index theorem has been extended to the equivariant case, to families of
operators, and to foliations.

The generalization this paper is concerned with is the index theorem for
formal deformation quantizations originally proved by Fedosov \cite{fe:book}
and (independantly) Nest--Tsygan \cite{nt}, also known as the algebraic index
theorem. In principle this is an abstract theorem computing the pairing of
$K$-theory classes with the cyclic cocycle given by the unique trace on a
formal deformation quantization of a symplectic manifold. However, as shown
in \cite{nt96}, the nomenclature ``index theorem'' may be justified by the
fact that in the case of the cotangent bundle with its canonical symplectic
form and the deformation quantization given by the asymptotic
pseudo-differential calculus, one recovers the original Atiyah--Singer
index theorem.

Here, we prove an algebraic index theorem for formal deformation quantizations
of symplectic orbifolds and derive from it its analytic version, the
well-known Kawasaki index formula for  orbifolds. Let us explain and state
the theorem in some more detail.

In our setup, see Section \ref{Sec:preli} for more details, a symplectic
orbifold $X$ is modeled by a proper \'etale groupoid
$G:G_1\rightrightarrows G_0$ with an invariant, nondegenerate two-form on
$G_0$. Consequently, we consider formal deformation quantizations of the
convolution algebra of $G$. As in \cite{ta:quantization}, these are
constructed by a crossed product construction of a $G$-invariant formal
deformation $\calA^\hbar$ of $G_0$ by $G$, denoted $\calA^\hbar\rtimes G$.
This is the starting point for the algebraic index theorem.

By computing the cyclic cohomology of the algebra $\calA^\hbar\rtimes G$ we have
given a complete classification of all traces in our previous paper
\cite{nppt}: the dimension of the space of traces equals the number of
connected components of the so-called inertia orbifold $\tilde{X}$
associated to $X$. Therefore, in contrast to the case of symplectic manifolds,
there is no unique (normalized) trace on the deformation of the convolution
algebra of $G$. In this paper, we construct a particular trace $\Tr$
using the notion of a twisted trace density, the appropriate
generalization of a trace density to orbifolds.
Standard constructions in K-theory, see Section \ref{Sec:preli}, yield a map
\begin{displaymath}
   \Tr_* :K_{\text{\tiny \rm orb}}^0(X)\rightarrow\C((\hbar)),
\end{displaymath}
associated to the trace $\Tr$, called the index map. Here, $K_{\text{\tiny \rm orb}}^0(X)$ is the Grothendieck group generated by isomorphism classes of so called orbifold vector bundles.
The index theorem proved in this paper expresses the value of this
map on a virtual orbifold vector bundle $[E] -[F] $ in terms of the
characteristic classes of $E$, $F$, the orbifold $X$ and the chosen
deformation.
\\[2mm] \noindent
{\bf Theorem}
  {\it Let $G$ be a proper \'etale Lie groupoid representing a
  symplectic orbifold $X$. Let $E$ and $F$ be $G$-vector bundles
  which are isomorphic outside a compact subset of $X$.
  Then the following formula holds for the index of $[E]-[F]$:}
\begin{displaymath}
  \Tr_* ([E]-[F]) = \int_{\tilde X} \frac1m
  \frac{\Ch_\theta \big(\frac{R^E}{2\pi i}-
  \frac{R^F}{2\pi i}\big)}{
  \det\big(1-\theta^{-1}\exp \big(-\frac{R^\perp}{2\pi i}\big)\big)}
  \hat{A}\Big( \frac{R^T}{2\pi i} \Big) \,
  \exp\Big(-\frac{\iota^*\Omega}{2\pi i \hbar}\Big).
\end{displaymath}
  \vspace{2mm}

In this formula, the right hand side is a purely topological
expression that we now briefly explain: Let $B_0$ be the space of
objects of the groupoid modeling the inertia orbifold $\tilde{X}$,
and $\iota:B_0\rightarrow G_0$ the canonical map explained in
Section \ref{Sec:preli}. This groupoid has a canonical cyclic
structure $\theta$ acting on the fibers the vector bundles $E$ and
$F$ as well as the normal bundle to the map $\iota$. This enables
one to define the twisted Chern character $\Ch_\theta$, see
Section \ref{Sec:locrr}. The coefficient function $\frac 1m$ is
defined in terms of the order of the local isotropy groups (see
Sections \ref{Sec:intorb} and \ref{SubSec:TwVB}). The factor $\hat{A}(\ldots)$ is the
standard characteristic class of the symplectic manifold $B_0$
associated to the bundle of symplectic frames, and the factor
$\det(\ldots)$ is the inverse of the twisted Chern character
associated to the symplectic frame bundle on the normal bundle to
$\iota:B_0\rightarrow G_0$. Finally, $\Omega$ is the
characteristic class of the deformation quantization of $G_0$.

Let us make several remarks about this Theorem. First, indeed notice that
the right hand side is a formal Laurent series in $\C((\hbar))$. Only when the
characteristic class of the deformation quantization is
trivial, it is independent of $\hbar$. This happens for example in the case
of a cotangent bundle $T^*X$ to a reduced orbifold $X$, with the canonical
deformation quantizatioin constructed from asymptotic pseudodifferential
calculus, in which case we show in Section \ref{Sec:algind}
that the left hand side equals the index of an elliptic operator on $X$.
This is exactly the Kawasaki index theorem \cite{Kaw:IEOVM}, originally
derived from the Atiyah--Segal--Singer $G$-index theorem.
An alternative proof using operator algebraic methods has been given by
Farsi \cite{farsi}.

Notice that, since we work exclusively with the convolution algebra of the
group\-oid and its deformation, the abstract theorem above also holds for
nonreduced orbifolds. In the reduced case it proves a conjecture of
Fedosov--Schulze--Tarkhanov \cite{fst}. Although they work with the
different algebra of invariants of $\calA^\hbar(G_0)$ instead of the crossed
product, one can show that in the reduced case the two are Morita equivalent,
allowing for a precise translation.

The methods used to prove our main result are related to the proof of the
algebraic index theorem on a symplectic manifold in \cite{FFS}. One of the
main tools in that paper is the construction of a trace density for a
deformation quantization on the underlying symplectic manifold. We generalize
the construction of such a trace density to symplectic orbifolds. Hereby,
our approach is inspired by the localization behavior of cyclic cocycles
on the deformed algebra,
which has been discovered in our previous paper \cite{nppt}, and includes
an essential new feature, a local ``twisting'' on the inertia groupoid.

Our paper is set up as follows. After introducing some preliminary material
in Section \ref{Sec:preli}, we introduce in Setion \ref{Sec:twtr} the
concept of a twisted trace density. In Section \ref{Sec:twhoch}, we
construct certain local twisted Hochschild cocyles, which in Section \ref{Sec:cotwtrd}
will be ``glued" to a twisted trace density. Thus, we obtain a trace
on the deformation quantization of the groupoid algebra, and consequently the index map on $K_{\text{\tiny \rm orb}}^0(X)$. In Section \ref{Sec:locrr}, we use
Chern-Weil theory on Lie algebras to determine  the cohomology class
of the trace density defined in Section \ref{Sec:cotwtrd}, and obtain a
local algebraic Riemann-Roch formula for symplectic orbifolds. Finally,
in Section \ref{Sec:algind} we prove the above theorem and use it to give
an algebraic proof of the classical Kawasaki index theorem for orbifolds.
In the Appendix, we provide some material needed for the proof our index
theorem. More precisely, in Appendix \ref{Sec:appa} we determine
several Lie algebra cohomologies and in Appendix \ref{Sec:appb} we explain the
asymptotic pseudodifferential calculus and its relation to
deformation quantization. \\[2mm]
\noindent{{\bf Acknowledgement:}} M.P.~and H.P.~acknowledge
financial support by the Deutsche Forschungsgemeinschaft. H.P.~and
X.T.~would like to thank the Fachbereich Mathe\-matik of
Goethe-Universit\"at at Frankfurt/Main for the hospitality during
their visits. X.T.~would like to thank Dmitry Fuchs, Ilya Shapiro,
and Alan Weinstein for helpful suggestions. H.P.~is supported by
EC contract MRTN-CT-2003-505078 (LieGrits).

%% file: preli.tex
\section{Preliminaries}
\label{Sec:preli}
In this section we briefly recall some of the essential points of \cite{nppt}.
For any orbifold $X$, we choose a proper \'etale Lie groupoid
$G_1\rightrightarrows G_0$ for which $G_0/G_1\cong X$. Such a groupoid always
exists and is unique up to Morita equivalence, although being \'etale is not
Morita invariant. The groupoid is used to describe the differential geometry
of the orbifold. For example, an orbifold vector bundle is equivalent to a
$G$-vector bundle, that is, a vector bundle $E\rightarrow G_0$ with an
isomorphism $s^*E\cong t^*E$. The notion of a $G$-sheaf is similarly defined.

The Burghelea space of $G$, also called ``space of loops'', is defined by
\begin{displaymath}
   B_0=\{g\in G_1 \mid s(g)=t(g)\},
\end{displaymath}
where $s,t:G_1\rightarrow G_0$ are the source
and target map of the groupoid $G$. Denote by $\iota$ the canonical embedding
$B_0\hookrightarrow G_1$. The groupoid $G$ acts on $B_0$ by
conjugating the loops and this defines the associated inertia groupoid
$\Lambda G:=B_0 \rtimes G$. This groupoid turns out also to be proper and
\'etale, and therefore models an orbifold, $\tilde{X}$, called the inertia
orbifold.

An important property of the inertia groupoid which will be essential in this
paper is the existence of a \textit{cyclic structure} \cite{crainic}: there
is a canonical section $\theta:\Lambda G_0\rightarrow\Lambda G_1$ of both
the source and target maps of $\Lambda G$, given by $g\mapsto \theta_g$.
Here $\theta_g=g$, viewed as a morphism from $g$ to $g$.

The convolution algebra of $G$ is defined as the vector
space $\calC^\infty_{\tiny\rm c}(G_1)$ with the following product:
\begin{equation}
\label{convolution}
  (f_1*f_2)(g)=\sum_{g_1g_2=g}f_1(g_1)f_2(g_2),\quad
  f_1,f_2 \in \calC^\infty_{\tiny\rm c}(G_1), ~ g\in G_1.
\end{equation}
The convolution algebra  is denoted by $\calA \rtimes G$, where $\calA$
here and everywhere in this article means the $G$-sheaf of smooth
functions on $G_0$.

When the orbifold $X$ is symplectic, one can choose $G$ in such a way that
$G_0$ carries a symplectic form $\omega$ which is invariant, i.e.,
$s^*\omega=t^*\omega$. For such a groupoid, we choose a $G$-invariant
deformation quantization of $G_0$, giving rise to a $G$-sheaf of algebras
$\calA^\hbar$ on $G$. The crossed product algebra $\calA^\hbar\rtimes G$
defines a formal deformation quantization of the convolution
algebra $\calA \rtimes G$ with its canonical noncommutative Poisson
structure induced by $\omega$. Recall that $\calA^\hbar\rtimes G$
is defined as the vector space $\Gamma_\text{\tiny\rm c}(G_1,s^*\calA^\hbar)$
with product
\begin{equation}
\label{crossedprd}
  [a_1 \star_c a_2 ]_g= \sum_{g_1 \, g_2=g}\big( [a_1]_{g_1}g_2 \big)
  [a_2]_{g_2},
  \quad
  a_1,a_2 \in \Gamma_\text{\tiny\rm c}(G_1,s^*\calA^\hbar), ~ g\in G,
\end{equation}
where $[a]_g$ denotes the germ of $a$ at $g$. Closely related is the algebra
$\Gamma_\text{\tiny\rm inv,c}(\calA^\hbar)$ of $G$-invariant sections, first
considered in \cite{P1}. In fact, when $X$ is reduced, the two are Morita
invariant, cf \cite[Prop. 6.5.]{nppt}, the equivalence bimodule being given
by $\calA^\hbar_c(G_0)$.

Let us finally mention that for every $G$-sheaf $\calS$
on $G_0$, we will denote the sheaf $s^{-1} \calS =t^{-1} \calS$ also by
$\calS$. Since $G$ is assumed to be proper \'etale, this will be convenient
and not lead to any misunderstandings.

\subsection{Traces on $\calA^\hbar\rtimes G$}
The starting point of the present article is the classification of traces
on the deformed convolution algebra $\calA^\hbar\rtimes G$ in \cite{nppt},
which we now briefly recall. A trace on the algebra $\calA^\hbar\rtimes G$
is an $\hbar$-adically continuous linear functional
$\Tr:\calA^\hbar\rtimes G\rightarrow\K$, where $\K$ denotes the field
of Laurent series $\C((\hbar))$, such that
\begin{displaymath}
  \Tr(a\star_c b)=\Tr(b\star_c a).
\end{displaymath}
Of course, a trace on an algebra is nothing but a cyclic cocycle of degree $0$.
Moreover, the space of traces on $\calA^\hbar\rtimes G$ is in bijective
correspondence with the space of traces on the extended deformed
convolution algebra $\calA^{((\hbar))}\rtimes G$, where
$\calA^{((\hbar))} = \calA^\hbar \otimes_{\C [[\hbar]]} \K$.
One of the main results of \cite{nppt} now asserts that
\begin{equation}
  HC^p(\calA^{((\hbar))}\rtimes G)\cong
  \bigoplus_{l\geq 0} H^{p-2l}\left(\tilde{X},\K\right),
\end{equation}
and therefore $HC^0$ equals $H^0(\tilde{X},\C)\otimes\K$. From
this it follows that the number of linear independent traces on
$\calA^\hbar\rtimes G$ equals the number of connected components
of the inertia orbifold $\tilde{X}$. In \cite{nppt}, a
construction of all these traces was given using a \v{C}ech-like
description of cyclic cohomology, however the resulting formulas
are not easily applicable to index theory. Therefore, we start in
Section \ref{Sec:twtr} by giving an alternative, more local,
construction of traces on deformed groupoid algebras using the
notion of a twisted trace density.
\subsection{The index map}
\label{Sec:indexmap}
Here we briefly explain the construction of the index map, given a trace
$\Tr:\calA^\hbar\rtimes G\rightarrow\K$ on the deformed convolution algebra.
As is well known,
a trace $\tau : A\rightarrow \Bbbk$ on an algebra $A$ over a field $\Bbbk$
induces a map in $K$-theory $\tau_*:K_0(A)\rightarrow\Bbbk$ by taking the
trace of idempotents in $M_n(A)$. In our case, we obtain a map
\begin{displaymath}
   \Tr_* :K_0(\calA^{((\hbar))}\rtimes G)\rightarrow\C((\hbar)).
\end{displaymath}
The inclusion
$\calA^\hbar\rtimes G\hookrightarrow \calA^{((\hbar))}\rtimes G$
induces a map
$K_0 (\calA^\hbar\rtimes G)\rightarrow K_0(\calA^{((\hbar))}\rtimes G)$.
Since $\calA^\hbar\rtimes G$ is a deformation quantization of
the convolution algebra of $G$, one has, by rigidity of $K$-theory
\begin{displaymath}
  K_0(\calA^\hbar\rtimes G)\cong K_0(\calA\rtimes G)\cong K^0(G).
\end{displaymath}
Here, $K^0(G)$ is the Grothendieck group of isomorphism
classes of $G$-vector bundles on $G_0$, also called the orbifold
$K$-theory $K_\text{\tiny\rm orb}^0(X)$. To be precise, the isomorphism
$K_0(\calA\rtimes G)\cong K^0(G)$, supposing that $X$ is compact, associates
to any $G$-vector bundle $E\rightarrow G_0$ the projective
$\calC^\infty_\text{\rm c}(G)$-module $\Gamma_\text{\rm c}(E)$, where
$f\in \calC^\infty_\text{c}(G)$ acts on $s\in\Gamma_\text{\rm c}(E)$ by
\begin{displaymath}
  (f\cdot s)(x)=\sum_{t(g)=x} f(g)\, s(x),~\text{for $x\in G_0$}.
\end{displaymath}
Putting all these maps together, the trace $\Tr$ defines a map
\begin{equation}
\label{indexmap}
  \Tr_*: K_{\text{\tiny \rm orb}}^0(X)\rightarrow\C((\hbar)).
\end{equation}
As we have seen in the previous paragraph, traces on
$\calA^{((\hbar))}\rtimes G$ are highly non-unique on general orbifolds, and
all induce maps as in \eqref{indexmap}. The trace we will construct in this
paper using the index density has the desirable property that it has support
on each of the components of $\tilde{X}$. Consequently, we will refer to the
map \eqref{indexmap} induced by this trace as the index map and the main
theorem proved in this paper gives a cohomological formula for the value on
a given  orbifold vector bundle $E$. Notice that because of the support
properties of this ``canonical trace'', this theorem in principle solves
the index problem for any other element in $HC^0(\calA^{((\hbar))}\rtimes G)$.
\subsection{Integration on orbifolds}
\label{Sec:intorb}
Let $\Gamma$ be a finite group acting by diffeomorphims on a
smooth manifold $M$ of dimension $n$. Consider the orbifold $X = M/\Gamma$,
and assume $X$ to be connected.
Recall that $X$ carries a natural stratification by orbit types
(see \cite[Sec.~4.3]{P3}). Denote by $X^\circ$ the principal stratum
of $X$ and by $M^\circ \subset M$ its preimage under the canonical projection
$M \rightarrow X$.
Now let $\mu $ be an $n$-form on $X$ or in other words a $\Gamma$-invariant
$n$-form on $M$. The integral $\int_X \mu$ then is defined by
\begin{equation}
\label{Eq:intorb}
  \int_X \mu := \frac{m}{|\Gamma|} \int_M \mu,
\end{equation}
where $m\in \N$ is the order of the isotropy group $\Gamma_x$ of some
point $x \in M^\circ$. Observe that by the slice theorem $m$ does not
depend on the particular choice of  $x$, since $X$ is connected.
Let us mention that formula \eqref{Eq:intorb} is motivated by the fact that
on the one hand the covering $M^\circ \rightarrow X^\circ$ has exactly
$\frac{|\Gamma|}{m}$ sheets and on the other hand the singular set
$X\setminus X^\circ$ has measure $0$.

Assume now to be given an arbitrary orbifold $X$. To define an integral over
$X$ choose a covering of $X$ by orbifold charts and a subordinate smooth
partition of unity. Locally, the integral is defined by the above
formula for orbit spaces. These local integrals are glued together globally
by the chosen partition of unity. The details of this construction are
straightforward.

%% file: twtr.tex
\section{Twisted trace densities}
\label{Sec:twtr}
Since $\calA^\hbar$ is a $G$-sheaf of algebras, the pull-back sheaf
$\iota^{-1}\calA^\hbar$ is canonically a $\Lambda G$-sheaf. This implies
that every stalk $\iota^{-1}\calA^\hbar_g,~g\in B_0$ has a canonical
automorphism given by the action of the cyclic structure
$\theta_g\in{\rm Aut}(\iota^{-1}\calA^\hbar_g)$. Alternatively, $\theta$
defines a canonical section of the sheaf
$\underline{\rm Aut}(\iota^{-1}\calA^\hbar)$ of automorphisms of
$\iota^{-1}\calA^\hbar$. This facilitates the following definition.
\begin{definition}
A $\theta$-twisted trace density on a $\Lambda G$-sheaf of algebras
$\iota^{-1}\calA^\hbar$ is a sheaf morphism
$\psi:\iota^{-1}\calA^\hbar\rightarrow\Omega^{\mbox{\tiny top}}_{\Lambda G}$
which satisfies
\begin{equation}
\label{twisteddensity}
  \psi(ab)-\psi(\theta(b)a)\in d\Omega^{\mbox{\tiny top}-1}_{\Lambda G}.
\end{equation}
\end{definition}
In this definition $\Omega^{\mbox{\tiny top}}_{\Lambda G}$ means the sheaf
of top degree de Rham forms on every connected component of $B_0$. Therefore,
the integral of a twisted trace density over $B_0$ is well defined. Notice
that although $\Omega^\bullet_{\Lambda G}$ is a $\Lambda G$-sheaf,
the action of $\theta$ is trivial, and therefore the twisting in the
definition only involves $\iota^{-1}\calA^\hbar$. As we will see later this
involves the normal bundle to the embedding $\iota :B_0\hookrightarrow G_1$.
\begin{proposition}
\label{traced}
When $\psi$ is a $\theta$-twisted trace density, the formula
\begin{displaymath}
  \Tr(a)=\int_{B_0}\psi(a_{|B_{0}}), \quad a\in\calA^{((\hbar))}\rtimes G,
\end{displaymath}
defines a trace on the deformed convolution algebra
$\calA^{((\hbar))}\rtimes G$.
\end{proposition}
\begin{proof}
The map $a\mapsto a_{|B_0}$ is the degree $0$ part of a natural morphism
\begin{displaymath}
  C_\bullet (\calA^{((\hbar))}\rtimes G)\rightarrow
  \Gamma_\text{\tiny\rm c}\left(B_\bullet ,\sigma^{-1}
  \calA^{((\hbar))}_{G_0^{\bullet+1}}\right),
\end{displaymath}
of complexes called ``reduction to loops''. Using the notation from
\cite[Sec.~2.6]{nppt}), $\calA^{((\hbar))}_{G_0^p}$ denotes here the
sheaf $(\calA^{((\hbar))})^{\boxtimes p}$ on the cartesian product $G_0^p$,
$B_p$ the Burghelea space given by
\begin{displaymath}
  B_p =\{ (g_0,\cdots,g_p) \in G^{p+1} \mid t(g_0)=s(g_p), \cdots,
  t(g_p) = s(g_{p-1})\},
\end{displaymath}
and $\sigma : B_p \rightarrow G_0^{p+1}$
the map $(g_0,\cdots ,g_p) \mapsto (s(g_0),\cdots , s(g_p))$.
Moreover, the differential on
$C_\bullet$ is the standard Hochschild differential on
$\calA^{((\hbar))}\rtimes G$, and the differential on the right hand side
is given in \cite{crainic,nppt}. In the following, we will only need
the first one
\begin{equation}
\label{map}
  d_0-d_1:\Gamma_\text{\tiny\rm c}\left(B_1,\sigma^{-1}
  \calA^{((\hbar))}_{G_0^2}
  \right)\rightarrow \Gamma_\text{\tiny\rm c}
  \left(B_0,\iota^{-1}\calA_{G_0}^{((\hbar))} \right).
\end{equation}
At the level of germs it is given by
\begin{eqnarray*}
  d_0[a_0,a_1]_{(g_0,g_1)}&=&[a_0g_1a_1]_{g_0g_1}, \\
  d_1[a_0,a_1]_{(g_0,g_1)}&=&[a_1g_0a_0]_{g_1g_0}.
\end{eqnarray*}
The complex $(C_\bullet(A),b)$ calculates the Hochschild homology
$HH_\bullet (A)$ for any $H$-unital algebra $A$. Since $HH_0(A)=A/[A,A]$,
a trace on $A$ is nothing but a linear functional on $HH_0(A)$.

It was proved in \cite[Prop.~5.7]{nppt} that reduction to loops induces a
quasi-isomorphism of complexes, in particular commutes with the
differentials. Therefore, to construct a trace on $\calA^\hbar\rtimes G$,
it suffices to construct a linear functional on the vector space
$\Gamma_\text{\tiny\rm c}(B_{0},\iota^{-1}\calA_{G_0}^{((\hbar))})$
which vanishes on the image of the map \eqref{map}.
>From the definition of the simplicial
operators $d_0$ and $d_1$ above, we see that the germ at $g\in B_0$
of such a section is given by
\begin{displaymath}
  \sum_{g_0g_1=g}\left([a_0]g_1[a_1]-[a_1]g_1g^{-1}[a_0]\right)=
  \sum_{g_0g_1=g}\left([a_0]g_1[a_1]_{g_1}-\theta([a_1])_{g_1}[a_0])\right),
\end{displaymath}
where, to pass over to the right hand side,  we have used that $\calA^{((\hbar))}$ is a
$G$-sheaf of algebras. By the defining property of the twisted trace
density such elements will be mapped into
$d\Omega^{\text{\tiny \rm top}-1}_{B_0}$. Integrating, it follows from
Stokes' theorem that the functional given by integrating the trace density
vanishes on the image of $d_0-d_1$. Combined with the restriction to $B_0$,
it follows that the formula in the proposition defines a trace.
\end{proof}
\begin{remark}
A warning is in order here, as the formula in Proposition \ref{traced} for
the trace appears to suggest that the trace only depends on the restriction
of a formal power series $a\in \calC^\infty_\text{\tiny\rm c}(G)[[\hbar]]$
to $B_0$. This is not true, since the reduction to loops maps an element
$a$ to its \textit{germ} at $B_0$, viewed as an element of
$\Gamma_\text{\tiny\rm c}(B_0,\iota^{-1}\calA^\hbar)$, and
otherwise the twisting condition would be trivial. In fact, we will see
that the normal bundle to $B_0\subset G_1$ plays an essential role in
the computations below.
\end{remark}
Concluding, to construct a trace, it suffices to construct a twisted trace
density. Inspired by the recent construction in \cite{FFS} of a trace
density from a certain Hochschild cocycle, we aim for a similar construction
in a twisted Hochschild complex.

%% file: twhoch.tex
\section{A twisted Hochschild cocycle}
\label{Sec:twhoch}
\subsection{Twisted Hochschild cohomology}
\label{twhc}
Let $A$ be an algebra equipped with an automorphism $\gamma\in{\rm Aut}(A)$.
There is a standard way to twist the Hochschild homology and cohomology of
$A$ by $\gamma$: recall that the Hochschild homology and cohomology
$H_\bullet (A,M)$ and $H^\bullet(A,M)$ are defined with values in any
bimodule $M$. Then one defines
\begin{equation}
\label{def:twh}
  HH_\bullet^\gamma(A):=H_\bullet(A,A_\gamma),
  \quad HH^\bullet_\gamma(A):=H^\bullet(A,A^*_\gamma),
\end{equation}
where $A_\gamma$ is the $A$-bimodule given by $A$ with bimodule structure
\begin{displaymath}
  a_1\cdot a\cdot a_2=a_1a\gamma(a_2),~ a_1,a_2\in A,~a\in A_\gamma,
\end{displaymath}
and $A_\gamma^*$ denotes its dual. Clearly this definition is functorial
with respect to automorphism preserving algebra homomorphisms. In the case
of cohomology, let us write out the standard complex computing this twisted
cohomology.

On the space of cochains
$C^p(A)={\rm Hom}(A^{\otimes p},A^*)\cong{\rm Hom}(A^{\otimes (p+1)},\K)$,
introduce the differential $b_\gamma:C^p\rightarrow C^{p+1}$ by
\begin{equation}
\label{twistedhochschild}
\begin{split}
  (b_\gamma f)(a_0\otimes\cdots \otimes a_{p+1})=
  \sum_{i=0}^p & (-1)^i  f(a_0,\cdots,a_ia_{i+1},\cdots,a_{p+1})\\
  &+(-1)^{p+1}f(\gamma(a_{p+1})a_0,a_1,\cdots,a_p).
\end{split}
\end{equation}
One checks that $b_\gamma^2=0$, and the cohomology of the resulting
cochain complex is called the twisted Hochschild cohomology of $A$.
When $\gamma=1$, this definition reduces to the ordinary Hochschild
cohomology $HH^\bullet(A)$.
\subsection{The local model}
\label{lm}
Let $V=\R^{2n}$ equipped with the standard symplectic form
$$\omega=\sum_{i=1}^ndp_i\wedge dq_i,$$ in coordinates
$(p_1,\ldots,p_n,q_1,\ldots,q_n)\in\R^{2n}$. The Weyl algebra
$\mathbb{W}_{2n}$ of $(\R^{2n},\omega)$ is given over the field
$\K = \C ((\hbar))$ by the vector space
$\K[p_1,\cdots p_n,q_1,\cdots,q_n]$ of polynomials in
$(p_1,\cdots,p_n,q_1,\cdots,q_n)$ with the Moyal product defined by
\begin{displaymath}
  a\star b= m\big(\exp\hbar\alpha (a\otimes b)\big).
\end{displaymath}
Here, $m:\mathbb{W}_{2n}\otimes\mathbb{W}_{2n}\rightarrow \mathbb{W}_{2n}$
is the commutative product on polynomials and
\begin{displaymath}
  \alpha(a\otimes b)=\sum_{i=1}^n \frac{\partial a}{\partial p_i}\otimes
  \frac{\partial b}{\partial q_i}-\frac{\partial b}{\partial p_i}\otimes
  \frac{\partial a}{\partial q_i}
\end{displaymath}
denotes the action of the Poisson tensor on
$\mathbb{W}_{2n}\otimes\mathbb{W}_{2n}$.
The Weyl algebra is in fact a functor from symplectic vector spaces to
unital algebras over $\K$, which implies that there is a canonical action
of $\Sp_{2n}$, the group of real linear symplectic transformations, on
$\mathbb{W}_{2n}$ by automorphisms. This induces an action of
 $\mathfrak{sp}_{2n}(\K):= \mathfrak{sp}_{2n} \otimes_\R \K$,
where  $\mathfrak{sp}_{2n}$ is the Lie algebra of $\Sp_{2n}$, by derivations.
In fact, this action is inner, $\mathfrak{sp}_{2n}(\K)$ being identified
with the degree two homogeneous polynomials in $\mathbb{W}_{2n}$ acting
by the commutator. In particular, a linear symplectic transformation
$\gamma\in\Sp_{2n}$ acts on the Weyl algebra by automorphisms and we can
consider the twisted Hochschild homology of $\mathbb{W}_{2n}$.
\begin{proposition}
\label{computationtwh}
{\rm (}Cf.~\cite{AFLS:HIAWAGF}{\rm )}
Let $\gamma$ be a linear symplectomorphism of $\R^{2n}$, and $2k$ the dimension
of the fixed point space of $\gamma$. Then the twisted Hochschild homology
of $\W_{2n}$ is given by
\begin{equation}
  HH_p^\gamma(\mathbb{W}_{2n})=
  \begin{cases}
     \K& \text{for $p=2k$},\\
     0& \text{for $p\neq 2k$}.
  \end{cases}
\end{equation}
\end{proposition}
\begin{proof}
Although the proposition can be extracted from \cite{nppt}, let us give the
(standard) argument. First note that there is a decomposition of the symplectic
vector space $W=\R^{2n}$ as $W=W^\gamma\oplus W^\perp$, where
$W^\gamma  = \ker (1-\gamma)$ is the fixed point space of $\gamma$ and
$W^\perp  = \im (1-\gamma)$ its symplectic orthogonal.
Since $\gamma \in \Sp_{2n}$, this is a symplectic decomposition, and by
assumption $\dim(W^\gamma)=2k$. Choose a symplectic basis
$(y_1,\cdots,y_{2k})$ of $W^\gamma$ and extend it to a symplectic basis
$(y_1,\cdots,y_{2n})$ of $W$.
Observe now that the filtration by powers of $\hbar$ induces a spectral
sequence with $E^0$-term the classical twisted Hochschild homology of the
polynomial algebra $A_{2n}=\K[p_1,\cdots,p_n,q_1, \cdots ,q_n]$. To compute
this homology, one chooses a Koszul resolution of $A$:
\begin{displaymath}
  0\longleftarrow A\stackrel{m}{\longleftarrow}A^e
  \stackrel{\partial}{\longleftarrow}\ldots
  \stackrel{\partial}{\longleftarrow} A^e \otimes \Lambda^{2n-1}W^*
  \stackrel{\partial}{\longleftarrow} A^e \otimes \Lambda^{2n}W^*
  \longleftarrow 0,
\end{displaymath}
where the differential $\partial$ is defined in the usual way by
\begin{displaymath}
\begin{split}
\partial ( a_1\otimes a_2 \otimes &\, dy_{i_1} \wedge\ldots\wedge dy_{i_p})=\\
  &\sum_{j=1}^p \, (-1)^j \, (y_{i_j}a_1\otimes a_2-a_1\otimes y_{i_j}a_2)
  \otimes dy_{i_1}\wedge\ldots\wedge\widehat{dy}_{i_j}\wedge\ldots\wedge
  dy_{i_p}.
\end{split}
\end{displaymath}
This is a projective resolution of $A$ in the category of $A$-bimodules.
Since by definition, cf. \eqref{def:twh},
$HH_\bullet^\gamma(A)={\rm Tor}^{A^e}_\bullet(A_\gamma,A)$, we apply the
functor $A_\gamma \otimes_{A^e} -$ to this resolution and compute the
cohomology. This yields the complex
\begin{displaymath}
  0\longleftarrow A_\gamma\stackrel{m}{\longleftarrow} A_\gamma \otimes W^*
  \stackrel{\partial}{\longleftarrow}\ldots \stackrel{\partial}{\longleftarrow}
  A_\gamma \otimes \Lambda^{2n-1}W^*  \stackrel{\partial}{\longleftarrow}
  A_\gamma \otimes \Lambda^{2n}W^* \longleftarrow 0
\end{displaymath}
with differential
\begin{displaymath}
  \partial ( a\otimes dy_{i_1}\wedge\ldots\wedge dy_{i_p}) =
  \sum_{j=1}^p \, (-1)^j \, (y_{i_j} a-\gamma( y_{i_j})a)\otimes
  dy_{i_1}\wedge\ldots\wedge\widehat{dy}_{i_j}\wedge\ldots\wedge dy_{i_p}.
\end{displaymath}
Since $y_1,\cdots y_{2k}$ form a basis of $V^\gamma$, their contribution
in the differential will vanish. The decomposition
$W=W^\gamma\oplus W^\perp$ yields a decomposition of exterior products
\begin{displaymath}
  \Lambda^lW^*=\bigoplus_{p+q=l}\Lambda^p(W^\gamma)^*
  \oplus\Lambda^q(W^\perp)^*,
\end{displaymath}
and from the previous remark we see that the differential on the
$\Lambda^p(W^\gamma)^*$-part is zero. In the $q$-direction one finds a
direct sum of degree shifted Koszul complexes of the ring $A$ associated
to the regular sequence
$(y_{k+1}-\gamma(y_{k+1}),\cdots,y_{2n}-\gamma(y_{2n}))$.
Recall that a Koszul complex of a regular sequence
$\mathbf{x}=(x_1,\cdots,x_n)$ in an algebra $A$  has homology concentrated
in degree $0$ equal to $A/I$, where $I=(\mathbf{x})A$. In our case this
ideal is exactly the vanishing ideal of $W^\gamma$, and one finds
\begin{displaymath}
  HH^\gamma_\bullet (A)=\Omega^\bullet_{A_{2k}/\K},
\end{displaymath}
where $A_{2k}=\K[y_1,\cdots,y_{2k}]$ and $\Omega^\bullet_{A_{2k}/\K}$ is the
algebra of K\"{a}hler differentials. Notice that this argument proves
a twisted analogue of the Hochschild--Kostant--Rosenberg theorem.

We therefore find $E^0_{p,q}=\Omega^{p+q}_{A_{2k}/\K}$, and the differential
$d^0:E^0_{p,q}\rightarrow E^0_{p,q-1}$ is the algebraic version of
Brylinski's Poisson differential \cite{byi:poisson} on the symplectic
variety $W^\gamma$. Using the symplectic duality transform this complex
is isomorphic to the algebraic de Rham complex of $W^\gamma$ shifted by
degree. Therefore one finds
\begin{displaymath}
  E^1_{p,q}=\K \quad \text{iff $p+q=2k$}.
\end{displaymath}
The spectral sequence collapses at this point and the result follows.
\end{proof}
It is not difficult to check that the cocycle
\begin{displaymath}
  c_{2k}=\sum_{\sigma\in \operatorname{S}_{2k}}\operatorname{sgn}(\sigma)
  1\otimes y_{\sigma(1)}\otimes\ldots\otimes y_{\sigma(2k)}
\end{displaymath}
is a generator of $HH^\gamma_{2k}(\mathbb{W}_{2n})$.
Notice that the cocycle only involves the basis elements on $W^\gamma$.

\subsection{The external product}
Let $A$ and $B$ be algebras equipped with automorphisms $\gamma_A$
and $\gamma_B$. Since, as we have seen, the twisted Hochschild
cohomology is nothing but the ordinary Hochschild cohomology with
values in a twisted bimodule, the external product in Hochschild
cohomology (see e.g.~\cite[Sec. 9.4.]{weibel}) yields a map
\begin{equation}
\label{extprod}
  \#:HH^p_{\gamma_A}(A)\otimes HH^q_{\gamma_B}(B)
  \rightarrow HH^{p+q}_{\gamma_A\otimes\gamma_B}(A\otimes B).
\end{equation}
Let us describe its construction. Denote by $B_\bullet(A)$ and $B_\bullet(B)$
the bar resolution of $A$ resp.~$B$ in the category of bimodules.
Their tensor product $B_\bullet(A)\otimes B_\bullet(B)$ carries the structure
of a bisimplicial vector space, which, by the Eilenberg--Zilber theorem is
chain homotopy equivalent to its diagonal. But the diagonal
$\operatorname{Diag} (B_\bullet(A)\otimes B_\bullet(B))$ is naturally
isomorphic to the Bar complex of $A\otimes B$. Explicitly, the
Eilenberg--Zilber theorem is induced by the so called Alexander--Whitney map,
see \cite[Sec.8.5.]{weibel}, which in our case gives rise to maps
$f_{pq}:B_{p+q}(A\otimes B)\rightarrow B_p(A)\otimes B_q(B)$. Taking the
$\Hom$ in the category of bimodules over $A\otimes B$ to
$A_\gamma\otimes B_\gamma$, combined with the natural map
\begin{displaymath}
  \Hom_{A^e}(B_p(A),A_{\gamma_A})\otimes {\rm Hom}_{B^e}(B_p(B),B_{\gamma_B})
  \rightarrow {\rm Hom}_{A^e\otimes B^e}(B_p(A\otimes B),A_{\gamma_A}
  \otimes B_{\gamma_B})
\end{displaymath}
yields a map $C^p(A)\otimes C^q(B)\rightarrow C^{p+q}(A\otimes B)$ commuting
with the twisted differentials, which induces \eqref{extprod}. The explicit
expression on the level of Hochschild cocycles can be read off from the
Alexander--Whitney mapping. We will only need the following special case:

Assume $\gamma_A=1$, and for notational simplicity put $\gamma:=\gamma_B$.
As one easily observes from the definition of the differential
\eqref{twistedhochschild}, a cocycle of degree $0$ in the Hochschild
cochain complex is nothing but a $\gamma$-twisted trace. Recall that a
$\gamma$-twisted trace on a $\K$-algebra $B$ with a fixed automorphism
$\gamma\in{\rm Aut}(B)$ is a linear functional
$\tr_\gamma:B\rightarrow \K$ such that
\begin{equation}
\label{twtrace}
  \tr_\gamma(b_1 b_2)=\tr_\gamma(\gamma(b_2)b_1).
\end{equation}
Taking the external product with such a cocycle yields the following:
\begin{lemma}
\label{ext}
Let $\tau_{k}$ be a Hochschild cocycle on an algebra $A$ of degree $k$. Then,
for a $\gamma$-twisted trace $\tr_\gamma$ on a $\K$-algebra  $B$ with an
automorphism $\gamma$, the formula
\begin{equation*}
  \tau_{k}^\gamma \big( (a_0 \otimes b_0) \otimes \ldots \otimes
  (a_{k} \otimes b_{k}) \big) :=
  \tau_{k} (a_0 \otimes \ldots\otimes a_{k}) \,
  \operatorname{tr}_\gamma (b_0 \cdots b_{k}),
\end{equation*}
defines a $\gamma$-twisted Hochschild cocycle of degree $k$ on $A\otimes B$.
\end{lemma}
\begin{proof}
Of course, one can prove this by checking that
$\tau^\gamma_k=\tau_k\#\tr_\gamma$, a fact which can be read off from the
precise form of the Alexander--Whitney map, and using that the map $\#$
passes to cohomology, as in \eqref{extprod}. However it can also be done
by a direct computation:
\begin{displaymath}
\begin{split}
  (b_\gamma\tau^\gamma_{k})&
  \big((a_0\otimes b_0)\otimes\cdots\otimes(a_{k+1}\otimes
  b_{k+1})\big)= \\
  & = \sum_{i=0}^{k} \, (-1)^i\, \tau_{k}(a_0\otimes\cdots\otimes
  a_i a_{i+1} \otimes\cdots\otimes a_{k+1}) \tr_\gamma (b_0\cdots b_{k+1}) \\
  & \hspace{2em} + (-1)^{k+1}\, \tau_ {k}(a_{k+1}
  a_0\otimes\cdots\otimes a_{k}) \tr_\gamma(\gamma(b_{k+1})b_0\cdots b_k)\\
  &=(b\tau_k)(a_0\otimes\cdots\otimes a_{k+1})\tr_\gamma(b_0\cdots b_{k+1}) =0,
\end{split}
\end{displaymath}
since $\tau_k$ is a Hochschild cocycle. Here, $b$ denotes the ordinary
(untwisted) Hochschild coboundary operator on the complex $C^k(A)$ and we
have used the twisted trace property \eqref{twtrace} of $\tr_\gamma$.
\end{proof}

\subsection{The cocycle}
In this subsection we will construct a cocycle of degree $2k$ in the
twisted Hochschild complex of $\mathbb{W}_{2n}$, where the twisting is
induced by $\gamma\in \Sp_{2n}$. We use the notation from the proof of
Prop.~\ref{computationtwh} and recall, in particular, the symplectic
decomposition $\C^n=\C^k \oplus \C^{(n-k)}$, where $\C^k=\ker(1-\gamma)$
and $\C^{(n-k)}={\rm Im}(1-\gamma)$. By this, the Weyl algebra has
form $\mathbb{W}_{2n}=\mathbb{W}_{2k}\otimes \mathbb{W}_{2(n-k)}$, and we
can use the external product to construct the Hochschild cocycle.
For notational simplicity, put
$\mathbb{W}_{2n}^\text{\tiny\rm T}:=\mathbb{W}_{2k}$ and
$\mathbb{W}_{2n}^\perp:=\mathbb{W}_{2(n-k)}$. The
computation of the twisted Hochschild cohomology of
$\mathbb{W},~\mathbb{W}^\text{\tiny\rm T}$ and $\mathbb{W}^\perp$ follows from
Prop.~\ref{computationtwh}. The twisted Hochschild cohomology of
$\mathbb{W}_{2n}$ is one dimensional and concentrated in degree $2k$.
Since the twisting induced by $\gamma$ is trivial on
$\mathbb{W}^\text{\tiny\rm T}_{2n}$,
its (ordinary, i.e., untwisted) Hochschild cohomology is also one dimensional
and in degree $2k$. Finally, $\mathbb{W}^\perp_{2n}$ has cohomology
concentrated in degree $0$ and equal to $\K$. Therefore, the only way to
construct a nontrivial twisted Hochschild cocycle on $\mathbb{W}_{2n}$ is
by taking the external product of an untwisted Hochschild cocycle of degree
$2k$ on $\mathbb{W}_{2n}^\text{\tiny\rm T}$ with a twisted trace on
$\mathbb{W}^\perp_{2n}$.

In \cite{FFS}, a Hochschild cocycle of degree $2k$ on $\mathbb{W}_{2k}$ was
constructed. Let us recall its definition. For $0\leq i\neq j\leq 2k$, denote
by $\alpha_{ij}$ the Poisson tensor on $\C^{2k}$ on the $i$'th and $j$'th
slot of the tensor product $\mathbb{W}_{2k}^{\otimes(2k+1)}$:
\begin{displaymath}
\begin{split}
   \alpha_{ij}(a_0\otimes\cdots \otimes a_{2k})=
   \frac{1}{2}\sum_{l=1}^k&\left(a_0\otimes\cdots \otimes
   \frac{\partial a_i}{\partial p_l}\otimes\cdots\otimes
   \frac{\partial a_j}{\partial q_l}\otimes\cdots\otimes
   a_{2k}\right.\\ &\left. -a_0\otimes\cdots\otimes
   \frac{\partial a_i}{\partial q_l}\otimes\cdots\otimes
   \frac{\partial a_j}{\partial p_l}\otimes\cdots\otimes a_{2k}\right).
\end{split}
\end{displaymath}
The operator $\pi_{2k}\in {\rm End}(\mathbb{W}_{2k}^{\otimes(2k+1)})$ is
defined as
\begin{displaymath}
  \pi_{2k}(a_0\otimes\cdots \otimes a_{2k})=
  \sum_{\sigma\in S_{2k}}\operatorname{sgn}(\sigma)\, a_0\otimes
  \frac{\partial a_1}{\partial y_{\sigma(1)}}\otimes\cdots\otimes
  \frac{\partial a_{2k}}{\partial y_{\sigma(2k)}}.
\end{displaymath}
Let $\mu_{2k}:\mathbb{W}_{2k}^{\otimes 2k+1}\rightarrow\K$ be the operator
$\mu_{2k}(a_0\otimes\cdots \otimes a_{2k})=a_0(0)\cdots a_{2k}(0)$, where
$a_i(0)$ is the constant term of $a_i$. With these operators at hand, define
\begin{equation}
  \tau_{2k} (a) =
  \mu_{2k} \int_{\Delta^{2k}} \prod_{0\leq i < j \leq 2k}
  e^{\hbar (2u_i -2u_j +1) \alpha_{ij}}
  \pi_{2k} (a) du_1 \wedge \ldots \wedge du_{2k} ,
\end{equation}
where
$a := a_0 \otimes a_1 \otimes \ldots \otimes a_{2k}\in
\mathbb W_{2k}^{\otimes 2k+1}$, and $\Delta^{2k}$ is the standard simplex
in $\R^{2k+1}$. As proved in \cite[Sec.~2]{FFS}, this defines a
nontrivial Hochschild cocycle of degree $2k$. It is an explicit cocycle
representative of the only non-vanishing cohomology class.

In the ``transverse direction'', i.e., on $\mathbb{W}^\perp_{2n}$, we need
a twisted trace. Fortunately, such traces have been constructed in
\cite{fe:g-index}. For this, we choose a $\gamma$-invariant complex
structure on $V^\perp$, identifying $V^\perp\cong\C^{n-k}$ so that
$\gamma\in U(n-k)$. The inverse Caley transform
\begin{displaymath}
  c(\gamma)=\frac{1-\gamma}{1+\gamma}
\end{displaymath}
is an anti-hermitian matrix, i.e., $c(\gamma)^*=-c(\gamma)$. With this,
define
\begin{displaymath}
  \tr_\gamma(a):=\mu_{2(n-k)}\left(\det^{-1}(1-\gamma^{-1})
  \exp\left(\hbar \, c(\gamma^{-1})^{ij}\frac{\partial}{\partial z^i}
  \frac{\partial}{\partial \bar{z}^j}\right)a\right),
\end{displaymath}
where $c(\gamma^{-1})^{ij}$ is the inverse matrix of $c(\gamma^{-1})$ and
we sum over the repeated indices $i,j=1,\ldots,n$. It is proved in
\cite[Thm.~1.1]{fe:g-index}, see also \cite[Lem.~7.3]{fst}, that this
functional is a $\gamma$-twisted trace density, i.e., satisfies equation
\eqref{twtrace}. Clearly, $\tr_\gamma(1)={\rm det}^{-1}(1-\gamma^{-1})$,
so we immediately see from Proposition \ref{computationtwh} that
$\tr_\gamma$ is independent of the choice of a complex structure.
This is also explicitly proved in \cite{fe:g-index}, but we view it as a
``cohomological rigidity''.

With the Hochschild cocycle $\tau_{2k}$ on $\mathbb{W}_{2k}$ and the
twisted trace density $\tau_\gamma$ on $\mathbb{W}_{2(n-k)}$ we can now
define the twisted Hochschild cocycle $\tau^\gamma_{2k}$ on
$\mathbb{W}_{2n}=\mathbb{W}_{2k}\otimes \mathbb{W}_{2(n-k)}$ of degree
$2k$ using the formula given in Lemma \ref{ext}.

As we have seen in section \ref{lm}, the Lie algebra $\mathfrak{sp}_{2n}$
acts on $\mathbb{W}_{2n}$ by derivations. Since $\gamma\in \Sp_{2n}$, it
will act on $\mathfrak{sp}_{2n}$ by the adjoint action. In the following,
we will be interested in the $\gamma$-fixed Lie subalgebra under this
action:
\begin{displaymath}
  \h:=\mathfrak{sp}_{2n}^\gamma\cong\mathfrak{sp}_{2k}\oplus
  \mathfrak{sp}_{2(n-k)}^\gamma.
\end{displaymath}
The isomorphism follows from the decomposition $\C^n=\C^k\oplus\C^{n-k}$,
which is a decomposition of representations of the cyclic group generated
by $\gamma$, $\C^k$ being the isotypical summand of the trivial one.
In general, $\h$ is a semisimple Lie subalgebra of $\mathfrak{sp}_{2n}$.
The following is the twisted analog of Theorem 2.2 in \cite{FFS} and lists
the properties of the cocycle obtained by means of Lemma \ref{ext}:
\begin{proposition}
\label{proptau}
  The cochain $\tau_{2k}^\gamma$ is a cocycle of degree $2k$ in the twisted,
  normalized Hochschild complex which has the following properties:
\begin{itemize}
\item[$i)$]
  The cochain $\tau_{2k}^\gamma$ is $\mathfrak{h}$-invariant which means that
  \begin{displaymath}
    \sum_{i=0}^{2k}\tau^\gamma_{2k}(a_0\otimes\ldots\otimes [a,a_i]
    \otimes\ldots\otimes a_{2k})=0 \quad \text{for all $a\in\mathfrak{h}$}.
  \end{displaymath}
\item[$ii)$]
   The relation
   $\tau^\gamma_{2k}(c_{2k})=\det^{-1}(1-\gamma^{-1})$ holds true.
\item[$iii)$]
   For every $a\in\mathfrak{h}$ one has
   \begin{displaymath}
     \sum_{i=1}^{2k}(-1)^i\tau^\gamma_{2k}(a_0\otimes\ldots\otimes
     a_{i-1}\otimes a\otimes a_i\otimes\ldots\otimes a_{2k-1})=0.
   \end{displaymath}
\end{itemize}
\end{proposition}
\begin{proof}
  It follows from Lemma \ref{ext} that $\tau_{2k}^\gamma$ is a twisted
  Hochschild cocycle of degree $2k$. Since $\tau_{2k}$ is normalized on
  $\mathbb{W}_{2n}$, see \cite[Thm.~2.2]{FFS}, the twisted cocycle
  $\tau^\gamma_{2k}$ is normalized as well. For the first property, we
  write an element in $a\in \h$ acting on $\mathbb{W}_{2n}$ as
  $a=x\otimes 1+1\otimes y$ with $x\in\mathfrak{sp}_{2k}$ and
  $y\in\mathfrak{sp}^\gamma_{2(n-k)}$. Consequently, we verify the
  equality in $i)$ for $x$ and $y$ seperately. For $x\in\mathfrak{sp}_{2k}$
  this is nothing but \cite[Thm.~2.2 $i)$]{FFS}. For
  $y\in \mathfrak{sp}^\gamma_{2(n-k)}$, one has
  \begin{displaymath}
    \sum_{i=0}^{2k}\tr_\gamma(a_0\cdots[y,a_i]\cdots a_{2k})=0,
  \end{displaymath}
  because $y$ commutes with $\gamma$. Property $ii)$ follows at once from
  the fact that $\tau_{2k}(c_{2k})=1$, cf.~\cite[Thm.~2.2 $ii)$]{FFS}.
  Finally, $ii)$ again splits into two parts for $a=x\otimes 1+1\otimes y$
  as above. The first part vanishes because of \cite[Thm.~2.2 $iii)$]{FFS}.
  The second is zero since $\tau_{2k}$ is normalized.
\end{proof}

%% file: cotwtrd.tex
\section{Construction of a twisted trace density}
\label{Sec:cotwtrd}
\subsection{Twisted Lie algebra cohomology}
\label{tlac} Let $\h\subset\g$ be an inclusion of Lie algebras.
Recall that for a $\g$ module $M$, the Lie algebra cochain complex
is given by $C^k(\g; M):={\rm Hom} (\Lambda^k\g,M)$, with
differential $\partial_{\mbox{\tiny Lie}}:C^k(\g;M)\rightarrow
C^{k+1}(\g;M)$ defined as
\begin{displaymath}
\begin{split}
  (\partial_{\mbox{\tiny Lie}}f)(x_1\wedge\ldots\wedge x_{k+1})&=
  \sum_{i=1}^{k+1}(-1)^{i+1}x_i\cdot f(x_1\wedge\ldots\wedge
  \hat{x}_i\wedge\ldots\wedge x_{k+1})\\ &+\sum_{i<j}(-1)^{i+j}f([x_i,x_j]
  \wedge\ldots\wedge\hat{x}_i\wedge\ldots\wedge\hat{x}_j\wedge\ldots\wedge
  x_{k+1}).
\end{split}
\end{displaymath}
This forms a complex, i.e., $\partial^2_{\text{\tiny Lie}}=0$, and
its cohomology is the Lie algebra cohomology $H^\bullet(\g;M)$.
Likewise, the relative Lie algebra cochain complex $C^k(\g,\h;
M)={\rm Hom}(\Lambda^k(\g/\h),M)^\h$ is the subcomplex of
$C^\bullet(\g;M)$ consisting of $\h$-invariant cochains vanishing
when any of the arguments is in $\h$. Its cohomology
$H^\bullet(\g,\h; M)$ is the relative Lie algebra cohomology with
coefficients in $M$.

For an algebra $A$ over $\mathbb{K}$, we denote by
$\mathfrak{gl}_N(A)$ the Lie algebra of $N\times N$ matrices with
entries in $A$. On equals footing, an $A$-bimodule $M$ yields a
$\mathfrak{gl}_N(A)$-module $\mathfrak{M}_N(M)$ of $N\times N$
matrices with entries in $M$ and module structure given by the
matrix commutator combined with the left and right $A$-module
structure. For $M=A_\gamma^*$, we therefore obtain the Lie algebra
complex $C^\bullet(\mathfrak{gl}_N(A);\mathfrak{M}_N(A^*_\gamma))$.
This is the twisted Lie algebra cohomology complex, and we denote
the differential by $\partial^\gamma_{\mbox{\tiny Lie}}$. Consider
now the map $\phi_N:C^\bullet(A,M)\rightarrow
C^\bullet(\mathfrak{gl}_N(A); \mathfrak{M}_N(A_\gamma^*))$ given by
\begin{displaymath}
  \phi_N(\tau)(x_1\wedge\ldots\wedge x_k)(x_0):=
  \sum_{\sigma\in S_{k}}{\rm sgn}(\sigma)\tau
  \left(\tr(x_0)\otimes\tr(x_{\sigma(1)})\otimes\ldots\otimes
  \tr(x_{\sigma(k)})\right),
\end{displaymath}
for $x_0\in\mathfrak{M}_N(A_\gamma)$, $x_1,\ldots,x_k\in\mathfrak{gl}_N(A)$,
and we used the matrix trace to define a map
$\tr:\mathfrak{gl}_N(A)\rightarrow A$. It is immediately clear by inspection
of the differentials that this defines a morphism
\begin{displaymath}
  \phi_N:\left(C^\bullet(A),b_\gamma\right)\rightarrow
  \left(C^\bullet(\mathfrak{gl}_N(A);\mathfrak{M}_N(A_\gamma^*)),
  \partial^\gamma_{\mbox{\tiny Lie}}\right)
\end{displaymath}
of cochain complexes. Using this morphism we define the following $2k$-cocycle
in the Lie algebra complex:
\begin{equation}
\label{liealgebracocycle}
  \Theta_{2k}^{N,\gamma}:=\phi_N(\tau_{2k}^\gamma).
\end{equation}
For $N=1$ we simply write $\Theta^\gamma_{2k}$ for this cocycle.

\subsection{The construction}
\label{constr}
We now come to the actual construction, which is in fact just a twisted
version of the construction in \cite{FFS}. As is well-known, any
deformation quantization of a symplectic manifold is isomorphic to a
Fedosov deformation. On the symplectic manifold $G_0$ this implies that
there is a resolution
\begin{displaymath}
  0\rightarrow \calA^\hbar(G_0)\longrightarrow\Omega^0(G_0,\mathcal{W}_{2n})
  \overset{D}{\longrightarrow}\Omega^1(G_0,\mathcal{W}_{2n})
  \overset{D}{\longrightarrow}\ldots
\end{displaymath}
of the space of global sections of the sheaf $\calA^\hbar$, where
$\Omega^p (G_0,\mathcal{W})$ is the space of $p$-forms with values in the
Weyl algebra bundle $\mathcal{W}_{2n}$, and $D$ is a so-called Fedosov
connection.
Recall that $\mathcal{W}_{2n}$ is the bundle of algebras defined by
\begin{equation}
  \mathcal{W}_{2n}=F^{2n}_{\Sp}\times_{\Sp_{2n}}\mathbb{W}_{2n},
\end{equation}
where $F^{2n}_{\Sp}$ denotes the bundle of symplectic frames on the
tangent bundle $TG_0$. Combining the wedge product of forms with the algebra
structure in the fibers of $\mathcal{W}$ turns
$\Omega^\bullet (G_0,\mathcal{W})$ into a graded algebra, with product denoted
by $\bullet$. By definition, a Fedosov connection is a connection on the
Weyl algebra bundle which is a derivation with respect to this product, i.e.,
\begin{displaymath}
  D(\alpha\bullet\beta)=(D\alpha)\bullet\beta+(-1)^{\mbox{\tiny deg}
  (\alpha)}\alpha\bullet(D\beta).
\end{displaymath}
The above resolution identifies $\calA^\hbar(G_0)$ as the space of flat
sections in $\Omega^0(G_0,\mathcal{W})$ compatible with its algebra
structure. The Fedosov connection can be decomposed as
\begin{equation}
\label{spl}
  D=\nabla+[A,-],
\end{equation}
where $\nabla$ is a symplectic connection on $TG_0$, i.e., $\nabla\omega=0$,
$[-,-]$ is the commutator with respect to the product $\bullet$,
and $A\in\Omega^1(G_0,\mathbb{W}_{2n})$. Since $D^2=0$, the quantity
\begin{equation}
\label{curvature}
\nabla A+\frac{1}{2}[A,A]=\Omega
\end{equation}
must be central, i.e., is a $\C[[\hbar]]$-valued two form. Since $G$ is a
proper \'etale Lie groupoid, the Weyl algebra bundle $\mathcal{W}_{2n}$ is
automatically a $G$-bundle, i.e., carries an action of $G$. We choose $D$
to be an invariant Fedosov connection with respect to this action. Then
this construction actually yields a resolution of $\calA^\hbar$ in the
category of $G$-sheaves. The associated symplectic connection $\nabla$ and
$\mathcal{W}_{2n}$-valued 1-form $A$ are therefore $G$-invariant.
Let us recall now that pull-back by $s$ (or equivalently by $t$) extends
$\mathcal{W}_{2n}$ (resp.~$D$,  $\nabla$, $A$) in a natural way to a
bundle (resp.~to connections respß.~a $1$-form) defined over $G_1$.
For convenience, we will denote in the following the thus obtained objects
by the same symbol as their restrictions to $G_0$.

Consider now a sector $\calO$ of $G$, that is, a minimal $G$-invariant
component of $B_0$. Using the natural embedding
$\iota:\calO\hookrightarrow G_1$, we can pull back the bundle of Weyl
algebras on $G_1$ to $\iota^*\mathcal{W}_{2n}$. As a pull-back, this bundle
inherits a natural Fedosov connection $\iota^*D=\iota^*\nabla+\iota^*A$
defined by
\begin{equation}
\label{pb}
  (\iota^*D)(\iota^*\alpha)=\iota^*(D\alpha),
\end{equation}
 with Weyl curvature
\begin{equation}
\label{pbcurv}
(\iota^*\nabla)(\iota^*A)+\frac{1}{2}[\iota^*A,\iota^*A]=\iota^*\Omega.
\end{equation}
By definition \eqref{pb}, restriction to $\calO$ maps flat
sections of $\Omega^\bullet(G_0,\mathcal{W}_{2n})$ to flat
sections over $\calO$. Combined with the natural inclusion
$\calA^\hbar(G_0)\rightarrow \Omega^0(G_0,\mathcal{W}_{2n})$ above
as flat sections with respect to the connection $D$, we obtain a
natural morphism of sheaves $\iota^{-1}\calA^\hbar\rightarrow
\iota^{-1}\mathcal{W}_{2n}$ which we write on sections as
$a\mapsto \iota^*a$. Indeed notice that the map
$\iota^*:\iota^{-1}\calA^\hbar(\calO)\rightarrow
\Omega^0(\calO,\iota^*\mathcal{W}_{2n})$ thus defined, depends on
the germ at $\calO$ of a section of $\calA^\hbar$ on $G_0$, since
its definition uses an embedding as a flat section of $D$ on $G_0$
before restricting to $\calO$. Applying this construction on every
sector
 $\calO\subset B_0$, one obtains a map
$\iota^*:\iota^{-1}\calA^\hbar(B_0)\rightarrow
\Omega^0(B_0,\iota^*\mathcal{W}_{2n})$.

After these preparations, define
\begin{equation}
\label{density}
\psi_D(a):=\frac{(-1)^k}{(2k)!}\Theta^\gamma_{2k}(\iota^*A\wedge\ldots\wedge
\iota^*A)(\iota^*a).
\end{equation}
\begin{lemma}
\label{mor}
  $\psi_D$ is a well-defined morphism of sheaves
  $\psi_D:\iota^{-1}\calA^\hbar\rightarrow \Omega^{\text{\tiny\rm top}}_{B_0}
  [[\hbar]]$ which depends only on the Fedosov connection $D$.
\end{lemma}
\begin{proof}
We have already showed that the germ of $\psi_D(a)$ at $x\in B_0$ depends
only on the germ of $a$ at $x$ in $G_0$, i.e., the morphism is defined on
$\iota^{-1}\calA^\hbar$. Next, we have to show that it is independent of
the choice of $A$, for a given Fedosov connection $D$. In general, the
splitting \eqref{spl} is unique up to a $\mathfrak{sp}_{2n}$-valued $1$-form
on $G_0$. Since the Fedosov connection $D$ is $G$-invariant, the restriction
to $B_0$ of the difference between two choices of $A$ is given by an
$\h$-valued 1-form on $B_0$. Therefore, it follows from
Proposition \ref{proptau} $iii)$ that $\psi_D$ only depends on the choice
of $D$.
\end{proof}
\begin{proposition}
\label{prop:trace-density}
  $\psi_D$ defines a $\theta$-twisted trace density on
  $\iota^{-1}\calA^\hbar$.
\end{proposition}
\begin{proof}
The proof is the same as in \cite{FFS}, except for the twisting: Consider
$\varphi:\iota^{-1}\calA^\hbar\otimes  \iota^{-1}\calA^\hbar\rightarrow
\Omega^{\text{\tiny top}-1}_{B_0}[[\hbar]]$ defined by
\begin{displaymath}
  \varphi(a\otimes b):=\frac{(-1)^k}{(2k)!}\Theta^\gamma_{2k}
  (\iota^*A\wedge\ldots\wedge \iota^*A\wedge \iota^*a)(\iota^*b).
\end{displaymath}
In the definition of $\iota^*a$, we extended $a$ to a flat section in
$\Omega^0(G_0,\mathcal{W}_{2n})$. By \eqref{pb}, $\iota^*a$ therefore is
flat with respect to $\iota^*D$ which means that
$\iota^*\nabla(\iota^*a)+[\iota^*A,\iota^*a]=0.$ Now we compute
\begin{displaymath}
\begin{split}
  d\varphi(a,b) =\, &
  (2k-1)\Theta^\gamma_{2k}((\iota^*\nabla)\iota^*A\wedge
  \ldots\wedge \iota^*A\wedge \iota^*a)(\iota^*b) \\
  &+\Theta^\gamma_{2k}(\iota^*A\wedge\ldots\wedge \iota^*A\wedge
  \iota^*\nabla \iota^*a)(b)\\
  &+\Theta^\gamma_{2k}(\iota^*A\wedge\ldots\wedge \iota^*A\wedge
  \iota^*a)(\iota^*\nabla \iota^*b)\\
  =\, &-\frac{2k-1}{2}\Theta^\gamma_{2k}([\iota^*A,\iota^*A]\wedge
  \iota^*A\ldots\wedge \iota^*A\wedge \iota^*a)(\iota^*b)\\
  &-\Theta^\gamma_{2k}(\iota^*A\wedge\ldots\wedge \iota^*A\wedge
  [\iota^*A,\iota^*a])(\iota^*b)\\
  &-\Theta^\gamma_{2k}(\iota^*A\wedge\ldots\wedge \iota^*A\wedge
  \iota^*a)([\iota^*A,\iota^*b])\\
  =\, &\Theta^\gamma_{2k}(\iota^*A\wedge\ldots\wedge
  \iota^*A)(\theta(\iota^*a)(\iota^*b)-
  (\iota^*b)(\iota^* a))\\
  =\, &\psi_D(\theta(a)b-ba),
\end{split}
\end{displaymath}
where we have used equation \eqref{pbcurv}. In the last step we have used
the fact that $\Theta^\gamma_{2k}$ is a twisted Lie algebra cocycle, together
with the fact that the connection $\iota^*A$ is $G$-invariant:
$\theta(\iota^*A)=\iota^*A$.
\end{proof}
Applying Proposition \ref{traced} we obtain a trace on the deformed
convolution algebra explicitly given by
\begin{equation}
\label{trace}
  \Tr(a)= \int_{B_0}\frac{1}{(2\pi i\hbar)^k} \psi_D(a)=
  \int_{B_0}\frac{1}{(2\pi i\hbar)^k} \Theta^\gamma_{2k}
  (\iota^*A\wedge\ldots\wedge \iota^*A)(\iota^*a).
\end{equation}
\begin{remark}
  In this formula and the remainder of this paper, $k$ will be regarded
  as an integer-valued function on $B_0$ (resp.~on $\tilde X$) which
  for every loop $g \in B_0$ gives half of the dimension of the fixed
  point space of the action of the isotropy group $G_x$ on $T_xG_0$, where
  $x = s(g)$. Clearly, $k$ is constant on each sector $\calO \subset B_0$
  by construction.
\end{remark}

\subsection{Twisting by vector bundles}
\label{SubSec:TwVB}
To evaluate the image of the general index map \eqref{indexmap}, it is
convenient to twist the construction of the density $\psi_D$ by an orbifold
vector bundle $E$. By definition, an orbifold vector bundle is a
$G$-vector bundle on $G_0$. The Fedosov construction can be twisted by a
vector bundle, see \cite{fe:book}, by tensoring the Weyl algebra bundle
$\mathcal{W}_{2n}$ with $\End (E)$. Choosing a connection $\nabla_E$ on
$E$, there is a flat Fedosov connection $D_E$ on
$\mathcal{W}_{2n}\otimes \End (E)$ which may be written as
$D_E=\nabla\otimes 1+1\otimes\nabla_E+[A_E,-]$, modifying \eqref{spl}.
The space of flat sections is a formal deformation $\calA^\hbar_E(G_0)$ of
the algebra $\Gamma(G_0,\End(E))$ of smooth sections of the bundle
$\End (E)$. In this section, we will construct traces on the crossed
product $\calA^\hbar_E\rtimes G$. Notice that since $\calA^\hbar_E\rtimes G$
is Morita equivalent to $\calA^\hbar\rtimes G$, its Hochschild and cyclic
homology are isomorphic and they have the same number of independent traces.

Since $G$ acts on the orbifold vector bundle, the restriction of $E$
(or more precisely $s^*E$) to $B_0$ carries a canonical fiberwise
action of the cyclic structure $\theta$. Let $V$ be the typical
fiber of $E$, which is a representation space of the finite group
$\Gamma$ generated by $\theta=\gamma$. In the local model, we
therefore switch from $\mathbb{W}_{2n}$ to the Weyl algebra with
twisted coefficients $\mathbb{W}_{2n}^V:=\mathbb{W}_{2n}\otimes \End
(V)$. As in Section \ref{tlac}, the $\mathbb{W}_{2n}^V$-bimodule
$\mathbb{W}^{V*}_{2n,\gamma}$, with the right action, in both
components, twisted by $\gamma$, yields a natural
$\gl(\mathbb{W}^V_{2n})$-module. Again, there is a natural morphism
$\phi_V:C^\bullet(\mathbb{W}^V_{2n};\mathbb{W}_{2n,\gamma}^{V*})
\rightarrow
C^\bullet(\gl(\mathbb{W}^V_{2n});\mathbb{W}^{V*}_{2n,\gamma})$ given
by
\begin{displaymath}
\begin{split}
  \phi_V(\tau)(M_1\otimes a_1,\ldots, M_k\otimes a_k)(M_0\otimes a_0):=&\\
  \sum_{\sigma\in S_{k}}{\rm sgn}(\sigma)\tau\left(a_0\otimes a_{\sigma(1)}
  \otimes\ldots \otimes a_{\sigma(k)}\right) &
  \tr_V\left(\gamma M_0 M_{\sigma(1)}\ldots M_{\sigma(k)}\right)
\end{split}
\end{displaymath}
for $M_0, M_1,\ldots M_k \in \End (V)$, $a_0,\ldots,a_k\in\mathbb{W}^V_{2n}$,
and $\tr_V$ is the canonical trace on $\End (V)$. The appearance of the
$\gamma$-twisted trace $\tr_{V,\gamma}$ can be explained by factorizing
this maps as the Morita equivalence
\begin{displaymath}
  \tr_{V,\gamma}:C^\bullet(\mathbb{W}_{2n};\mathbb{W}_{2n,\gamma}^{*})
  \rightarrow C^\bullet (\mathbb{W}^V_{2n};\mathbb{W}^{V*}_{2n,\gamma}),
\end{displaymath}
combined with the natural morphism
$C^\bullet(\mathbb{W}^V_{2n};\mathbb{W}^{V*}_{2n,\gamma})\rightarrow
C^\bullet(\gl(\mathbb{W}^V_{2n});\mathbb{W}^{V*}_{2n,\gamma})$ of
Section \ref{tlac}. For any $\Gamma$-representation $V$ we therefore
define
\begin{equation}
\label{Eq:theta}
  \Theta^{V,\gamma}_{2k}:=\phi_V(\tau^\gamma_{2k})\in C^{2k}
  (\gl(\mathbb{W}^V_{2n});\mathbb{W}^{V*}_{2n,\gamma}).
\end{equation}
Continuing as in Section \ref{constr}, the trace density
$\psi_{D_E}:\iota^{-1}\calA^\hbar_E\rightarrow\Omega^{\mbox{\tiny
top}}_{B_0}[[\hbar]]$ is defined by
\begin{displaymath}
  \psi_{D_E}(a)=\Theta^{V,\gamma}_{2k}
  (\iota^*A_E\wedge\ldots\wedge \iota^*A_E)(\iota^*a).
\end{displaymath}
The analogues of Lemma \ref{mor} and Proposition
\ref{prop:trace-density} are proved in exactly the same manner.
Therefore the integral over $B_0$ of this density gives a trace
$\Tr_E$ on $\calA^\hbar_E\rtimes G$. Its virtues lie in the
following proposition. Recall from (\ref{indexmap}) that the trace
defined in equation \eqref{trace} induces a map
$\Tr_*:K^0_{\text{\tiny\rm orb}}(X)\rightarrow\K$.
\begin{proposition}
\label{prop:k}
  Let $E$ and $F$ be orbifold vector bundles on $X$ which are isomorphic
  outside a compact subset.
  Then $[E]-[F]\in K^0_{\text{\tiny\rm orb}}(X)$, and one has
\begin{displaymath}
  \Tr_*([E]-[F])=  \int_{\tilde{X}} \frac{1}{(2\pi i \hbar)^k \, m}
  \big(\psi_{D_E}(1)-\psi_{D_F}(1)\big).
\end{displaymath}
Hereby, $m : \tilde X \rightarrow \N$ is the locally constant function which
coincides for every sector $\calO \subset B_0$ with $m_\calO$, the order
of the isotropy group of the principal stratum of $\calO/G\subset \tilde X$.
\end{proposition}
\begin{remark}
Note that both $\psi_{D_E}(1)$ and $\psi_{D_F}(1)$ are $G$-invariant, and
therefore define differential forms on $\tilde{X}$. Moreover, since $E$
and $F$ are isomorphic outside a compact subset of $X$,
$\psi_{D_E}(1)-\psi_{D_F}(1)$ is compactly supported on $\tilde{X}$.
\end{remark}
\begin{proof}
First recall from Section \ref{Sec:indexmap} that for every orbifold vector
bundle $E$, the section space $\Gamma_\text{c} (E)$ is a projective
$\calC_\text{c}^{\infty}(G)$-module. Hence, $\Gamma_\text{c}(E)$
defines a projection $e$ in the matrix algebra
$\mat_N(\calC_\text{c}^{\infty}(G))$ of $\calC_\text{c}^{\infty}(G)$, for
some large $N$. Observe that $e$ is a projection with respect to the
convolution product on $\calC_\text{c}^{\infty}(G)$, and generally is not a
projection valued matrix function.
According to \cite[Thm.~6.3.1]{fe:book}, the projection $e$ now has an
extension to a ``quantized'' projection $\hat{e}$ in
$\mat_N(\calA^{((\hbar))}\rtimes G)$.
By definition, $\Tr_*([E])$ equals $\Tr(\hat{e})$.

In the following, we will express $\Tr(\hat{e})$ in terms of an integral
over the inertia orbifold. We will do the computation locally, and use a
partition of unity later to glue the formulas.
When restricted to a small open set $U$ of the orbifold $X$, the
groupoid representing $U$ is Morita equivalent to a transformation
groupoid ${M}\rtimes \Gamma \rightrightarrows {M}$,
where ${M}$ is symplectomorphic to an open set of $\R^{2n}$
with the standard symplectic form, and $\Gamma$ is a finite
group acting by linear symplectomorphisms on ${M}$. The restriction
of $e$ to ${M}\rtimes \Gamma$ can be computed explicitly.

Any $\Gamma$-vector bundle $E$ on ${M}$ can be embedded
into a trivial $\Gamma$-vector bundle of high enough rank. Denote by $E'$ the
complement of $E$ in this trivial $\Gamma$-vector bundle, and let
$e^\Gamma$ be the projection valued $\Gamma$-invariant function corresponding
to the decomposition $E\oplus E'$. Let $\hat{e}^\Gamma$ be the
quantization of $e^\Gamma$. By \cite[Thm.~3]{cd:general-index},
one has $\psi_{D_E}(1)=\psi_{D_E+D_{E'}}(\hat{e}^\Gamma)$ for some
Fedosov connection $D_{E'}$  on $E'$.

Now let $\Pi$ denote the function
$\frac{1}{|\Gamma|}\sum_{\gamma\in\Gamma} \delta_{\gamma}$ on
${M}\rtimes \Gamma$, where $\delta_\gamma$ denotes the element
$(1,\gamma)$ in $\calC^{\infty}({M})\rtimes \Gamma$. Observe that $\Pi$ is
a projection with respect to the convolution product on
$\calC^{\infty}_\text{c}({M})\rtimes \Gamma$. The subspace
$M_\Pi:= \big(\calC_\text{c}^{\infty}({M})\rtimes\Gamma\big)\cdot\Pi$
generated by $\Pi$ in $\calC_\text{c}^{\infty}({M})\rtimes \Gamma$
is a left $\calC_\text{c}^{\infty}({M})\rtimes \Gamma$-module and a
right $\calC_\text{c}^{\infty}({M})^{\Gamma}$-module. Tensoring with
$M_\Pi$ thus defines a map from the projective modules on
$\calC_\text{c}^{\infty}({M})^{\Gamma}$ to
$\calC_\text{c}^{\infty}({M})\rtimes \Gamma$. Hence, it induces a map
on the corresponding projectors. It is not hard to see that
$e^\Gamma$ is mapped to $e$, and $e$ can be expressed by
$e^\Gamma$ as follows,
\begin{displaymath}
  e=\frac{1}{|\Gamma|}\sum_{\gamma\in \Gamma} e^{\Gamma} \delta_{\gamma}.
\end{displaymath}
The quantizations of $e$ and $e^\Gamma$ also satisfy this relation, hence
\begin{displaymath}
  \hat{e}=\frac{1}{|\Gamma|}\sum_{\gamma\in \Gamma}
  \hat{e}^{\Gamma} \delta_{\gamma}.
\end{displaymath}

The ``loop space" $B_0$ of the transformation groupoid
${M}\rtimes \Gamma$ is $\coprod_{\gamma\in \Gamma}
{M}^{\gamma}$, where ${M}^{\gamma}$ is the fixed point
set of $\gamma$. Hence, over $U$, the trace $\Tr^U_*([E])$ is computed
as follows,
\begin{displaymath}
\begin{split}
  \Tr^U_*([E])&= \int_{B_0}\frac{1}{(2\pi i \hbar)^k}
  \psi_D(\hat{e})=\sum_{\gamma\in \Gamma}\int_{{M}^{\gamma}}
  \frac{1}{(2\pi i \hbar)^k} \psi_D(\hat{e})\\
  &= \sum_{\gamma\in \Gamma} \int_{{M}^{\gamma}}
  \frac{1}{(2\pi i \hbar)^k \, |\Gamma|}\psi_D\big(\hat{e}^{\Gamma}\big)\\
  &=
  \sum_{\left<\gamma\right>}  \int_{{M}^{\gamma}}
  \frac{1}{(2\pi i \hbar)^k\, |\operatorname{C}(\gamma)|}
  \psi_D \big( \hat{e}^{\Gamma} \big)\\
  &=\sum_{\calO}\int_{\calO /\Gamma}\frac{1}{(2\pi i \hbar)^k \, m_\calO}
  \psi_D \big( \hat{e}^{\Gamma} \big).
\end{split}
\end{displaymath}
Hereby, $\left< \gamma \right>$ denotes the conjugacy class of $\gamma
\in \Gamma$, $\operatorname{C}(\gamma)$ is the centralizer of
$\gamma$, $\calO \subset B_0$ runs through the sectors of
$M\rtimes \Gamma$, and $\calO/\Gamma \subset \tilde{U}$ is the
quotient of $\calO$ in the inertia (sub)orbifold $\tilde{U}\subset
\tilde{X}$. Let us briefly justify the last equality in this
formula. To this end recall first the definition of an orbifold
integral from Section \ref{Sec:intorb} and second that each sector
$\calO/\Gamma$ coincides with the quotient of some connected
component of $M^\gamma$ by the centralizer $\operatorname{C}
(\gamma)$. Let $x$ be an element of the open stratum of this
connected component $M^\gamma_\alpha$ with respect to the
stratification by orbit types. By definition, $m_\calO$ is given
by the order of the isotropy group $\operatorname{C} (\gamma)_x$.
The orbit of  $\operatorname{C}(\gamma)$ through $x$ then has
$|\operatorname{C}(\gamma)|/m_\calO$ elements, and the integral
$\int_{\calO/\Gamma}\mu$ over an invariant form $\mu$ coincides
with
$\frac{m_\calO}{|\operatorname{C}(\gamma)|}\int_{M^\gamma_\alpha}
\mu$. This proves the claimed equality (cf.~also
\cite[Sec.~5]{fst}).

Using the fact that
$\psi_{D_E}(1)=\psi_{D_E+D_{E'}}(\hat{e}^{\Gamma})$, we thus obtain the
following equality:
\begin{equation}
\label{eq:loc-tr}
  \Tr^U_*([E])= \int_{\tilde{X}_U}\frac{1}{(2\pi i \hbar)^k \, m}\psi_{D_E}(1).
\end{equation}

Finally, note that the constructions of the projections $e$, and
$e^{\Gamma}$, and of the bimodule $M_\Pi$ are local with respect to $X$.
Therefore, one can glue together the local expressions (\ref{eq:loc-tr})
for the traces $\Tr^U_*$ by a partition of unity over $X$. This
completes the proof.
\end{proof}

%% file: locrr.tex
\section{A local Riemann--Roch theorem for orbifolds}
\label{Sec:locrr}
Our goal in this section is to use Chern-Weil theory in Lie
algebra cohomology to express the form $\psi_D(1)$ constructed in
Proposition \ref{prop:trace-density} by characteristic classes.

Recall that in Section 4, we have defined $\psi_D(a)$ as
$\Theta^{N,\gamma}_{2k}(\iota^*A\wedge\ldots\wedge
\iota^*A)(\iota^*a)$, where
\begin{displaymath}
   \Theta^{N,\gamma}_{2k}\in
   C^\bullet\big(\mathfrak{gl}_N(\mathbb{W}_{2n}),
   \mathfrak{gl}_N (\K) \oplus\mathfrak{sp}_{2k}(\K) \oplus
   \mathfrak{sp}_{2n-2k}^{\gamma}(\K) ; \mathfrak{M}_N(\W_{2n,\gamma}^*)
   \big)
\end{displaymath}
is a cocycle. Consider now the morphism $\ev_1:
\mathfrak{M}_N(\mathbb{W}_{2n,\gamma}^*)\rightarrow\mathbb{K}$,
which is the evaluation at the identity. Because of the nontrivial
$\gamma$-action, this is not a morphism of
$\mathfrak{gl}_N(\mathbb{W}_{2n})$-modules, but of
$\mathfrak{gl}_N(\mathbb{W}^\gamma_{2n})$-modules. Consequently, we
put
\begin{equation}
\label{liealgebras}
\begin{split}
  \g:=&\mathfrak{gl}_N(\mathbb{W}_{2n}^\gamma), \\
  \h:=& \mathfrak{gl}_N(\mathbb{K})\oplus \mathfrak{sp}_{2k}(\mathbb{K})\oplus
  \mathfrak{sp}_{2n-2k}^{\gamma}(\mathbb{K}),
\end{split}
\end{equation}
and consider the Lie algebra cohomology cochain complex
$C^{2k}(\g,\h;  \mathbb{K})$. It is in this cochain complex that we
explicitly identify the cocycle $\ev_1\Theta^{N,\gamma}_{2k}$.
Notice that this suffices to compute $\psi_D(1)$, since the
restricted connection $\iota^*A$ is by assumption invariant under
the action of the cyclic structure $\theta$.
\subsection{Chern-Weil theory}
In the following, we use Chern-Weil theory of Lie algebras to
determine the cohomology groups $H^p(\mathfrak{g},\h; \mathbb{K})$
for $p\leq 2k$.

First, recall the construction of Lie algebra Chern-Weil
homomorphism. As above, let $\frakg$ be a Lie algebra and $\frakh$ a
Lie subalgebra with an $\h$-invariant projection $\pr:\frakg\to
\frakh$. Define the curvature $C\in \Hom(\wedge^2 \frakg, \frakh)$
of $\pr$ by $C(u\wedge v):=[\pr(u),\pr(v)]-\pr([u,v])$. Let
$(S^\bullet\frakh^*)^\frakh$ be the algebra of $\frakh$-invariant
polynomials on $\frakh$ graded by polynomial degree. Define the
homomorphism $\chi:(S^\bullet\frakh^*)^\frakh\to
C^{2\bullet}(\frakg,\frakh;\K)$ by
\begin{displaymath}
  \chi(P)(v_1\wedge\cdots \wedge v_{2q})=
  \frac{1}{q!}\!\!\!\!\sum_{\sigma\in S_2q,\atop \
  \sigma(2i-1)<\sigma(2i)}\!\!\!\!(-1)^{\sigma}P\big( C(v_{\sigma(1)},\
  v_{\sigma(2)}), \cdots, C(v_{\sigma(2q-1)},\ v_{\sigma(2q)})\big).
\end{displaymath}
The right hand side of this equation defines a cocycle, and the
induced map in cohomology $\chi:(S^\bullet\frakh^*)^\frakh\to
H^{2\bullet}(\frakg,\frakh;\K)$ is independent of the choice of the
projection $\pr$. This is the Chern--Weil homomorphism.

In our case, that means with $\g$ and $\h$ given by Equation
\eqref{liealgebras},
the projection $\pr:\g\rightarrow\h$ is defined by
\begin{displaymath}
  \pr(M\otimes a):=\frac{1}{N}\tr(M)a_2+Ma_0,
\end{displaymath}
 where $a_j$ is the component of $a$ homogeneous of degree $j$ in $y$.
The essential point about the Chern--Weil homomorphism in this case is
contained in the following result.
\begin{proposition}
\label{prop:chern-weil}
  For $N\gg n$, the Chern-Weil homomorphism
  \begin{displaymath}
    \chi:(S^q\frakh^*)^{\frakh}\to H^{2q}(\frakg,\frakh;\K)
  \end{displaymath}
  is an isomorphism for $q\leq 2k$.
\end{proposition}
\begin{proof}
  The proof of this result goes along the same lines as the proof of
  Proposition 4.2 in \cite{FFS}, bearing in mind that $\h$ is semisimple.
  The only difference
  is that we need the following result on the cohomology
  $H^\bullet (\frakg;S^q\frakg^*)$ for $q >0$ proved in
  Corollary \ref{Cor:MainApp}:
\begin{displaymath}
 H^p (\mathfrak{gl}_N(\mathbb{W}_{2n}^{\gamma});
 S^q\mathfrak{M}_N(\mathbb{W}_{2n}^{\gamma}))=
 \begin{cases}
   0,& \text{for $j<2k$},\\
   \mathbb{K}^l, & \text{for $j=2k$}.
  \end{cases}
\end{displaymath}
\end{proof}
\subsection{The algebra $W_{k, n-k, N}^{\gamma}$}
By the isomorphism proved in Proposition \ref{prop:chern-weil}, it follows
that the cohomology class
\begin{displaymath}
 [\ev_1(\Theta_{2k}^{N, \gamma})]\in H^{2k}(\frakg , \frakh; \K )
\end{displaymath}
corresponds, under the Chern--Weil isomorphism, to a unique element
$P^{\gamma}_{2k}$ in ${(S^{2k}\frakh)^*}^{\frakh}$. To find this polynomial,
we restrict the Chern-Weil homomorphism to a Lie subalgebra
$W_{k, n-k, N}^{\gamma} \subset \g $.

Pick coordinates $(p_1,\ldots,p_k,q_1,\ldots q_k)$ of $\R^{2k}$ and
$(z_{k+1},\ldots z_{n-k})$ of $\C^{n-k}$ as in Section \ref{Sec:twhoch}.
Recall that $\gamma\in\Sp_{2n}$ acts trivially on the
$p_i,q_i,~1\leq i\leq k$. Moreover,  we assume that the
$z_j,~k+1\leq j\leq n$ are chosen to diagonalize $\gamma$. Let $K$ be the
polynomial algebra $K:=\mathbb{K}[q_1, \cdots, q_k, z_{k+1}, \cdots, z_n]$.
Introduce the Lie algebra $W_{k,n-k,N}:={\rm Der}(K)\rtimes\gl_N(K)$, the
semidirect product of polynomial vector fields by $\gl_N(K)$. This is the
Lie subalgebra of $\gl_N ( \mathbb{W}_{2n})$ of elements of the form
\begin{displaymath}
  \sum_i f_ip_i\otimes 1 +\sum_j g_j\bar{z}_j\otimes 1+\sum_k h_k\otimes M_k,
  \quad \text{with $f_i, g_j,h_k\in K$}.
\end{displaymath}
This Lie algebra has a natural $\gamma$-action, and
$W^\gamma_{k,n-k,N}$ is defined to be the $\gamma$-invariant part of
$W_{k,n-k,N}$. The Lie algebra
$\mathfrak{h}_1=W_{k, n-k,N}^{\gamma}\cap \mathfrak{h}$ is isomorphic to
\begin{displaymath}
  \mathfrak{gl}_k\oplus\mathfrak{gl}_{i_1}\oplus\cdots\oplus
  \mathfrak{gl}_{i_l}\oplus \mathfrak{gl}_N,
\end{displaymath}
where $\mathfrak{gl}_{i_s}$ corresponds to the eigenspace of the
$\gamma$-action on $\complex^{n}$ with a given eigenvalue.
Given this Lie subalgebra $\frakh_1\subset W_{n,n-k, N}^{\gamma}$, we
now consider the Chern-Weil homomorphism
\begin{equation}
\label{cwinv}
  \chi:(S^{q}\frakh^*_1)^{\frakh_1}\to H^{2q}\big(W_{k,
  n-k,N}^{\gamma},
  \frakh_1 ; \K \big).
\end{equation}

\begin{proposition}
\label{prop:w-inj}
  For $q\leq k$, the Chern-Weil homomorphism \eqref{cwinv}
  is injective.
\end{proposition}

\begin{proof}
By definition, $H^\bullet (W_{k, n-k, N}^{\gamma}, \frakh_1 ;
\mathbb{K})$ is the cohomology of the Cartan--Eilenberg cochain
complex
\begin{displaymath}
 \Big(\Hom\Big({\bigwedge}^\bullet
 \big( W_{k,n-k, N}^{\gamma}/\frakh_1
 \big), \mathbb{K}\Big) ^{\frakh_1},\partial_\text{\tiny \rm Lie}\Big).
\end{displaymath}
Write down a basis of $W_{k, n-k, N}$ as follows:
\begin{equation}
\label{eq:basis} q^{\alpha_i}z^{\beta_i}p_i\otimes 1,\
q^{\alpha_j}z^{\beta_j}\bar{z}_j\otimes 1,\
q^{\alpha_{st}}z^{\beta_{st}}\otimes E_{st},
\end{equation}
where $q^{\alpha_i}z^{\beta_i}p_i$, $q^{\alpha_j}z^{\beta_j}\bar{z}_j$, and
$q^{\alpha_{st}}z^{\beta_{st}}$ are polynomials. In the above formulas,
$\alpha_i, \alpha_j, \alpha_{st}, \beta_i, \beta_j, \beta_{st}$ are
multi-indices. If $\alpha_i=(\alpha_i^1, \cdots, \alpha_i^k)$, then
$q^{\alpha_i}:=q_1^{\alpha_i^1}\cdots q_k^{\alpha_i^k}$.
If $\beta_j=(\beta_j^1, \cdots, \beta_j^{n-k})$, then
$z^{\beta_j}:=z_1^{\beta_j^1}\cdots z_{n-k}^{\beta_j^{n-k}}$.
Finally, $E_{st},\ 1\leq s, t\leq N$ denotes the elementary matrix with $1$
at the $(s,t)$-position and $0$ everywhere else.
Since $\gamma$ acts diagonally on this basis, the $\gamma$-invariant
elements in \eqref{eq:basis} form a basis of $W^\gamma_{k,n-k,N}$.

Next, consider the $\frakh_1$-action on
$\Hom \big( {\bigwedge}^\bullet(W_{k, n-k, N}^{\gamma}),\mathbb{K}\big)$.
Note that the following elements act
diagonally on $W_{k, n-k, N}^{\gamma}$ and commute with each
other:
\begin{displaymath}
\sigma_1:=\sum_{i=1}^k q_ip_i\otimes 1,\ 1\leq i
\leq k \quad \text{and} \quad
\sigma_2:=\sum_{j=1}^{n-k}z_j\bar{z}_j\otimes 1,\ 1\leq j\leq n-k.
\end{displaymath}
Let us write down the formulas for the action of these elements.
\begin{enumerate}
\item $\sigma_1$ action.
\begin{enumerate}
\item
 $[\sigma_1, q^{\alpha_{i'}}z^{\beta_{i'}}p_{i'}\otimes 1]=
 \big(\big(\sum\limits_{l=1}^k \alpha_{i'}^l\big)-1\big)
 q^{\alpha_i}z^{\beta_i}p_i\otimes 1$,
 where $\alpha_{i'}^l$ is the $l$-th component of $\alpha_{i'}$.

\item
 $[\sigma_1, q^{\alpha_j}z^{\beta_j}\bar{z}_j\otimes 1]=
 \big( \sum\limits_{l=1}^k\alpha_j^{l}\big)
 q^{\alpha_j}z^{\beta_j}\bar{z}_j\otimes 1$,
 where $\alpha_j^l$ is the $l$-th component of $\alpha_j$.

\item
 $[\sigma_1, q^{\alpha_{st}}z^{\beta_{st}}\otimes E_{st}]=
 \big( \sum\limits_{l=1}^k\alpha_{st}^{l} \big)
 q^{\alpha_{st}}z^{\beta_{st}}\otimes E_{st}$, where
 $\alpha_{st}^l$ is the $l$-th component of $\alpha_{st}$.
\end{enumerate}
\item $\sigma_2$ action.
\begin{enumerate}
\item
 $[\sigma_2,  q^{\alpha_{i'}}z^{\beta_{i'}}p_{i'}\otimes 1]=
 \big(\sum\limits_{l=1}^{n-k}\beta_{i'}^l\big)
 q^{\alpha_{i'}}z^{\beta_{i'}}p_{i'}\otimes 1$,
 where $\beta_{i'}^l$ is the $l$-th component of $\beta_{i'}$.
\item
 $[\sigma_2, q^{\alpha_{j'}}z^{\beta_{j'}}\bar{z}_{j'}\otimes 1]=
 \big( \sum\limits_{l=1}^{n-k}\beta_{j'}^{l}-1\big)
 q^{\alpha_{j'}}z^{\beta_{j'}}\bar{z}_{j'}\otimes 1$,
 where $\beta_{j'}^l$ is the $l$-th component of $\beta_{j'}$.
 Note that $\bar{z}_i$ is not $\gamma$-invariant, and therefore does not
 belong to $W_{k, n-k, N}^{\gamma}$. If
 $q^{\alpha_{j'}}z^{\beta_{j'}}\bar{z}_{j'}\otimes 1$ is $\gamma$-invariant,
 $\sum\limits_{l=1}^{n-k}\beta_{j'}^{l}$ has to be greater than or equal to
 $1$.
\item $[\sigma_2, q^{\alpha_{st}}z^{\beta_{st}}\otimes E_{st}]=
 \big( \sum\limits_{l=1}^{n-k}\beta_{st}^{l} \big)
 q^{\alpha_{st}}z^{\beta_{st}}\otimes E_{st}$, where
 $\beta_{st}^l$ is the $l$-th component of $\alpha_{st}$.
\end{enumerate}
\end{enumerate}
In the following, we denote by $|\alpha_i|$  the sum of the
components of $\alpha_i$, and similarly for $|\beta_j|$.
Since $\sigma_1, \sigma_2$ are in $\frakh_1$, we know that only those
elements of
\begin{displaymath}
 \Hom \Big( {\bigwedge}^\bullet\big (W_{k, n-k, N}^{\gamma}/\frakh_1\big),
 \mathbb{K}\Big)
\end{displaymath}
which have eigenvalue $0$ will contribute to the $\frakh_1$-relative
cohomology.
Note that $\sigma_2$ acts only with nonnegative eigenvalues.
Therefore, we can reduce our considerations
to the $\sigma_2$-invariant Lie subalgebra of $W_{k, n-k, N}^{\gamma}$.
This Lie subalgebra will be denoted by $\tilde{W}_{k, n-k, N}^{\gamma}$
and has the following basis:
\begin{equation}
\label{eq:basis-inv}
 q^{\alpha_{i}}p_{i}\otimes 1,\
 q^{\alpha_{jj'}}z_j\bar{z}_{j'}\otimes 1,\ q^{\alpha_{st}}\otimes E_{st},
\end{equation}
where $z_j$ and $z_{j'}$ have the same eigenvalue for the $\gamma$-action,
and where $1\leq i\leq k$, $1\leq j,j'\leq n-k$, and $1\leq s,t\leq N$.
By passing to the quotient $\tilde{W}_{k, n-k, N}^{\gamma}/\frakh_1$,
no elements of the form
$q_{i'}p_i\otimes 1$, $z_j\bar{z}_{j'}\otimes 1$ or $1\otimes
E_{st}$ remain.
To compute the relative cohomology
$H^\bullet (W_{k, n-k, N}^{\gamma}, \frakh_1)$, we consider
the absolute cochain complex
\begin{displaymath}
  \Big( \Hom \big({\bigwedge}^\bullet(\tilde{W}_{k,n-k, N}^{\gamma}/\frakh_1),
  \K\big) ^{\frakh_1}, \partial_\text{\tiny \rm Lie}\Big).
\end{displaymath}
The Lie algebra $\tilde{W}_{k, n-k, N}^{\gamma}$ is closely related
to the Lie algebra $W_n$ considered in Theorem 2.2.4 of
\cite{fuchs:book}. In particular, by the same arguments from
invariant theory which show Lemma 1 in the proof of
\cite[Thm.~2.2.4]{fuchs:book} one concludes that
$\big(\Hom\big({\bigwedge}^\bullet (\tilde{W}_{k,n-k,
N}^{\gamma}/\frakh_1), \mathbb{K}\big)^{\frakh_1},
\partial_\text{\tiny \rm Lie}\big)$ is generated by even degree
polynomials. Let us explain this in the following in some detail.

First, denote the dual basis of (\ref{eq:basis-inv}) for
$\tilde{W}_{k,n-k, N}^{\gamma}$ as follows:
\begin{displaymath}
\begin{array}{ll}
 w_{\alpha_i}&=(q^{\alpha_{i}}p_{i}\otimes 1)^*\in
 \Hom(\tilde{W}_{k,n-k, N}^{\gamma}, \mathbb{K}),\\
 w_{\alpha_{jj'}}&=(q^{\alpha_{jj'}}z_j\bar{z}_{j'}\otimes
 1)^*\in \Hom(\tilde{W}_{k,n-k, N}^{\gamma}, \mathbb{K})\\
 w_{\alpha_{st}}&=(q^{\alpha_{st}}\otimes E_{st})^*\in
 \Hom(\tilde{W}_{k,n-k, N}^{\gamma}, \mathbb{K}).
\end{array}
\end{displaymath}
Moreover, if $I$ is a finite set of indices $\alpha_i$,
denote by $w_I$ the antisymmetric product
\begin{displaymath}
  w_I = w_{\alpha^1_{i_1}} \wedge \cdots  \wedge w_{\alpha^l_{i_l}},
\end{displaymath}
where $\alpha^1_{i_1} < \cdots < \alpha^l_{i_l}$ are the elements of
$I$ ordered by lexicographic (or some other fixed) order on $I$.
Likewise define $w_J$ and $w_S$ for each finite set of indices
$\alpha_{jj'}$ resp.~$\alpha_{st}$. Then, every element $\psi \in
\Hom({\bigwedge}^\bullet (\tilde{W}_{k,n-k, N}^{\gamma}/\frakh_1),
\K)^{\frakh_1}$ can be written as a linear combination
\begin{displaymath}
  \psi= \sum_{I,J,S} \psi_{I,J,S} \, w_I \wedge  w_J \wedge w_S ,
\end{displaymath}
where $\psi_{I,J,S} \in \K$,
$I$ runs through all finite sets of indices $\alpha_i$,
$J$ through the finite sets of indices $\alpha_{jj'}$ with
$|\alpha_{jj'}| \geq 1$ and $S$ through all finite sets of indices
$\alpha_{st}$ with $|\alpha_{st}|\geq 1$. Note that
the restriction on the elements of the
sets $J$ and $S$ comes from the fact that
$z_j\bar{z}_{j'}\otimes 1$ and $1\otimes E_{st}$ vanish in the quotient
$W_{k,n-k, N}^{\gamma}/\mathfrak{h}_1$.

Since $\psi$ is $\frakh_1$-invariant and the $\sigma_1$-action
is diagonal, the sum of the eigenvalues of the $\sigma_1$-action of
each single component of $w_I\wedge w_J\wedge w_S$
has to vanish, in case $\psi_{I,J,S} \neq 0$. In other words, if
$\psi_{I,J,S} \neq 0$, we have the following identity:
\begin{equation}
\label{eq:weight}
  \sum\limits_{\alpha_i \in I}(|\alpha_i|-1)+
  \sum\limits_{\alpha_{jj'} \in J}(|\alpha_{jj'}|)+
  \sum\limits_{\alpha_{st} \in S}(|\alpha_{st}|)=0.
\end{equation}
Suppose now that in $w_I\wedge w_J\wedge w_S$ there are $m_{I,J,S}$
elements of the form $w_i=(p_i\otimes1)^*$.  Then Equation (\ref{eq:weight})
entails
\begin{equation}
 \label{eq:sum-weight}
  m_{I,J,S} =\sum\limits_{
 \alpha_i>0\atop \alpha_i \in I}(|\alpha_i|-1)+
 \sum\limits_{\alpha_{jj'}>0\atop \alpha_{jj'} \in J}(|\alpha_{jj'}|)
 +\sum\limits_{\alpha_{st}>0\atop \alpha_{st} \in S}(|\alpha_{st}|).
\end{equation}
Furthermore, by the same arguments which show Lemma 1 in the proof
of \cite[Thm.~2.2.4]{fuchs:book}, one concludes from the
$\gl_k$-invariance of $\psi$, that the $w_i$'s have to be
paired with distinct terms of $w_{\alpha_k}$ or $w_{\alpha_{jj'}}$
or $w_{\alpha_{st}}$ such that $|\alpha_k|\geq2$, $|\alpha_{jj'}|,
|\alpha_{st}|\geq1$.
Thus, one has for $\psi_{I,J,S} \neq 0$ that $m_{I,J,S}$ is less than
the total number of terms of $w_{\alpha_k}$ , $w_{\alpha_{jj'}}$, and
$w_{\alpha_{st}}$  appearing in $w_I\wedge w_J\wedge w_S$ with
$|\alpha_k|$, $|\alpha_{jj'}|$, $|\alpha_{st}|\geq1$. Thus one has
\begin{equation}
\label{eq:pairing}
  m_{I,J,S} \leq\sum\limits_{
 \alpha_i>0\atop \alpha_i \in I}1+
  \sum\limits_{\alpha_{jj'}>0\atop \alpha_{jj'} \in J}1+\sum\limits_{\alpha_{st}>0\atop \alpha_{st} \in S}1.
\end{equation}

Equations (\ref{eq:sum-weight}) and (\ref{eq:pairing}) together
show that $w_{\alpha_i}$, $w_{\alpha_{jj'}}$, and
$w_{\alpha_{st}}$ cannot appear in $w_I\wedge w_J\wedge w_S$ for
nonvanishing $\psi_{I,J,S}$, if
$|\alpha_{i}|\geq 3$ or $|\alpha_{jj'}|\geq2$, or $|\alpha_{st}|\geq 2$.
Hence, by Equation (\ref{eq:sum-weight}) again, we have that $ m_{I,J,S}$
is equal to the number of terms $w_{\alpha_i}$, $w_{\alpha_{jj'}}$ and
$w_{\alpha_{st}}$ showing up in $w_I\wedge w_J\wedge w_S$ with
$|\alpha_i|=2$, $|\alpha_{jj'}|=1$, and $|\alpha_{st}|=1$. Therefore,
\begin{displaymath}
  m_{I,J,S}=\sum\limits_{|\alpha_i|=2}1+\sum\limits_{|\alpha_{jj'}|=1}1+
  \sum\limits_{|\alpha_{st}|=1}1.
\end{displaymath}
This shows that the degree of $w_I\wedge w_J\wedge w_S$  with
$\psi_{I,J,S} \neq 0$ is equal to $2m_{I,J,S}$, which is even. The
above arguments also show that
$\Hom\big({\bigwedge}^\bullet(\tilde{W}_{k,n-k,
N}^{\gamma}/\frakh_1), \K\big)^{\frakh_1}$ is generated by even
degree polynomials. This implies that the differential on the
cochain complex $\Hom\big({\bigwedge}^\bullet(\tilde{W}_{k,n-k,
N}^{\gamma}/\frakh_1), \K\big)^{\frakh_1}$ degenerates. Therefore,
in order to prove that $\chi$ is injective for $q\leq k$, it is
enough to show that for every nonzero polynomial $P \in
{(S^q\frakh_1)^{*}}^{\frakh_1}$ the image
\begin{displaymath}
  \chi (P) \in  H^{2q}(W_{k, n-k, N}^{\gamma},\frakh_1; \K)
   =
  \Hom\big({\bigwedge}^{2q}(\tilde{W}_{k,n-k, N}^{\gamma}/\frakh_1),\K\big)^{\frakh_1}
\end{displaymath}
does not vanish on ${\bigwedge}^{2q}(\tilde{W}_{k,n-k,
N}^{\gamma}/\frakh_1)$. But this follows from a straightforward
check.
\end{proof}

\subsection{Calculation of $\ev_1(\Theta_{2k}^{N, \gamma})$}
By Proposition \ref{prop:chern-weil} and Proposition \ref{prop:w-inj},
one obtains the following commutative diagram.
\begin{equation}
\label{diag:comm}
\begin{diagram}
  \node{(S^{k}\mathfrak{h}^*)^{\mathfrak{h}}}\arrow{e,t,T}{}\arrow{s,l,T}{\chi}
  \node{(S^k \mathfrak{h}_1^*)^{\mathfrak{h}_1}}\arrow{s,l,T}{\chi}\\
  \node{H^{2k}(\mathfrak{g}, \mathfrak{h};\K)}\arrow{e, t, T}{}
  \node{H^{2k}(W_{k,n-k, N}^{\gamma},\mathfrak{h}_1;\K)}
\end{diagram}
\end{equation}

The left vertical arrow in (\ref{diag:comm}) is the isomorphism
$\chi$ proved from Proposition \ref{prop:chern-weil}. The right
vertical arrow has been constructed in Proposition \ref{prop:w-inj}
and has been proved to be injective. The two horizontal arrows are
restriction maps. Since $\h$ and $\h_1$ have the same Cartan subalgebra
$\mathfrak{a}$ spanned by
\begin{displaymath}
  q_i p_i\otimes 1,~1\leq i\leq k, \quad z_j\bar{z}_j\otimes 1,~
  k+1\leq j\leq n, \quad  1\otimes E_r,~1\leq r\leq N,
\end{displaymath}
and since invariant polynomials are uniquely determined by their values on
$\mathfrak{a}$, the upper triangle is injective.
This implies that the lower horizontal
map is also injective. Therefore,  to determine a polynomial
$P_{k}^{\gamma}\in (S^{k}\frakh^*)^{\frakh}$ such that
$\chi(P_{k}^{\gamma})=\ev_1(\Theta_{2k}^{N, \gamma})$, one only needs
to work with the restriction of $\ev_1(\Theta_{2k}^{N, \gamma})$ to
$W_{k, n-k, N}^{\gamma}$.

Let $X=X_1\oplus X_2\oplus X_3\in \mathfrak{sp}_{2k}(\mathbb{K})
\oplus\mathfrak{sp}^\gamma_{2(n-k)}(\mathbb{K})\oplus
\mathfrak{gl}_N(\mathbb{K})=\mathfrak{h}$. Define
$(\hat{A}_{\hbar}\Ch_{\gamma}\Ch)_q\in (S^q\mathfrak{h}^*)^{\mathfrak{h}}$
to be the homogeneous terms of degree $q$ in the Taylor expansion of
\begin{displaymath}
  \left(\hat{A}_{\hbar}\Ch_{\gamma}\Ch\right)(X):=
  \hat{A}_{\hbar}(X_1)\Ch_{\gamma}(X_2)\Ch(X_3),
\end{displaymath}
where $\hat{A}_{\hbar}(X_1)$ is
$\det \big(\frac{{\tiny \hbar} X_1/2}{\sinh({\tiny \hbar} X_1/2)}
\big)^{\frac{1}{2}}$, and $\Ch_{\gamma}(X_2)$ is
$\tr_{\gamma}(\exp_\star (X_2))$ of the star exponential of $X_2$, and
$\Ch(X_3)$ is $\tr(\exp(X_3))$. Recall now that
$\ev_1(\Theta_{2k}^{N, \gamma})$ is given by the formula
\begin{displaymath}
  \ev_1\Theta_{2k}^{N,\gamma}(v_1, \cdots, v_{2k})=
  \Theta_{2k}^{N,\gamma}(v_1,\cdots, v_{2k})(1).
\end{displaymath}

\begin{theorem}
\label{thm:local-rr}
  For $N\gg n$, the following identity holds in $H^{2k}(\g,\h;\K)$:
\begin{displaymath}
  [\ev_1\Theta_{2k}^{N,\gamma}]=(-1)^k\chi((\hat{A}_{\hbar}\Ch\Ch_{\gamma})_k).
\end{displaymath}
\end{theorem}
\begin{proof}
Let us construct an invariant polynomial
$P^{\gamma}_{k}\in (S^\bullet\h^*)^\h$
as follows. Since $P^{\gamma}_{k}$ is required to be invariant under
the adjoint action, it is determined
by its value on the Cartan subalgebra $\mathfrak{a}$ which is spanned by
$p_iq_i\otimes 1, ~ 1\leq i\leq k$,
$z_j\bar{z}_j\otimes 1, ~ 1\leq j\leq n-k$, and
$1\otimes E_{l}, ~ 1\leq l\leq N$. Define $P_{k}^{\gamma}$ in
$(S^k\mathfrak{h}^*)^{\mathfrak{h}}$ to be the unique homogeneous
polynomial whose restriction to $\mathfrak{a}$ is
\begin{displaymath}
\begin{split}
  P_{k}^{\gamma} \, &
  (M_1\otimes a_1\otimes b_1, \cdots, M_k\otimes a_k\otimes b_k)
  =\tr(M_1\cdots M_k) \cdot \\
  & \cdot \mu_k\Big(\int_{[0,1]^k} \prod_{1\leq i \leq j\leq n}
  e^{\hbar \psi(u_i-u_j)\alpha_{ij}} (a_1\otimes\ldots\otimes a_k)
  du_1\cdots du_k\Big)\tr_{\gamma}(b_1\star\cdots \star b_k),
\end{split}
\end{displaymath}
with $\mu_k$ is as in Section 3.4.
\begin{remark}
If we define $P_{k}(M_1\otimes a_1,\cdots, M_k\otimes a_k)$ simply as
$P_{k}^{\gamma}(M_1\otimes a_1\otimes 1, \cdots, M_k\otimes
a_k\otimes 1)$, the polynomial $P_{k}$ is the same as the one in the proof of
Theorem 4.1 in \cite{FFS}. Clearly, $P_{k}^{\gamma}$ is the product of
$P_{k}$ and $\tr_{\gamma}$.
\end{remark}

In the following, we prove that
$[\ev_1\Theta_{2k}^{N, \gamma}]=(-1)^k\chi(P_{k}^{\gamma})$.
By $\mathfrak{h}$-invariance, it is enough to show the equality on
the Cartan subalgebra $\mathfrak{a}$. Let us define
$u_{ij}, v_{ir}, w_{is}\in W_{k, n-k, N}^{\gamma}$
for $i,j=1, \cdots, k$, $r=1, \cdots,N$, and $s=k+1,\cdots, n$
as follows:
\begin{displaymath}
\begin{array}{lll}
  u_{ij}:=
  \left\{\begin{array}{ll}\frac{1}{2}q_i^2p_i,&i=j,\\
          q_iq_jp_j& i\ne j,
         \end{array}\right.
  & v_{ir}:=q_i\otimes E_{r},&
  w_{is}:=q_iz_s\bar{z}_s.
\end{array}
\end{displaymath}
It is not difficult to check that
\begin{displaymath}
\begin{array}{lll}
  [p_i, u_{ij}]_{\hbar}=q_jp_j,& [p_i, v_{ir}]_{\hbar}=E_{r},&
  [p_i, w_{is}]_{\hbar}=z_s\bar{z}_s.
\end{array}
\end{displaymath}
Since $\pr(u_{ij})=\pr(v_{ir})=\pr(w_{is})=0$, one has
\begin{displaymath}
  C(p_i, u_{ij})=-p_iq_i, \quad C(p_i, v_{ir})=-E_{r}, \quad\text{and}\quad
  C(p_i, w_{is})= -z_s\bar{z}_s .
\end{displaymath}
Let $x_i$ be of the form $u_{jk}$ with $j\geq k$ or $v_{jr}$, or $w_{js}$.
Then it is straightforward to check that
\begin{displaymath}
  \ev_1\Theta_{2k}^{N,\gamma}(p_1\wedge x_1 \cdots \wedge p_k\wedge x_k)=
  (-1)^k P_{k}^{\gamma} \Big( \frac{x_1}{q_1}, \cdots, \frac{x_k}{q_k} \Big ).
\end{displaymath}
Note that in the definition of $P_{k}^{\gamma}$, we did not change the
component of $\tr_{\gamma}$, so the computation is the same as in the proof of
\cite[Thm.~4.1]{FFS}.

Now, we explicitly evaluate the polynomial $P_{k}^{\gamma}$ on the diagonal
matrices in the Cartan algebra. Let $X=Y+Z$, where
$Y:=\sum \nu_i q_ip_i +\sum \sigma_rE_{r}$ and $Z:=\sum \tau_s z_s\bar{z}_s$,
with $\nu_i,\sigma_r,\tau_s\in\mathbb{K}$.
Consider the generating function
\begin{displaymath}
  S(X):=\sum\limits_{m=1}^{\infty}\frac{1}{m!}P_{m}^{\gamma}(X,\cdots,X).
\end{displaymath}
Then
\begin{displaymath}
\begin{split}
  S(X)&=\sum_{m\geq 0}\frac{1}{m!}P_{m}^{\gamma}(X,\cdots, X) \\
  &=\sum_{m\geq 0}\frac{1}{m!}\sum_{0\leq l \leq m}\frac{m!}{l!(m-l)!}P_{l}
  (\underbrace{Y,\cdots,Y}_{l})\tr_{\gamma}
  (\underbrace{Z\star\cdots\star Z}_{m-l})\\
  &=\sum_{m\geq 0}\sum_{0\leq l \leq m}\frac{1}{l!(m-l)!}P_{l}
  (Y,\cdots,Y)\tr_{\gamma}(Z^{\star (m-l)})\\
  &=\sum_{l\geq 0}\frac{1}{l!}P_{l}(Y,\cdots,Y)\sum_{k=m-l\geq 0}
  \frac{1}{k!}\tr_{\gamma}(Z^{\star k}).
\end{split}
\end{displaymath}
According to \cite[Thm 4.1]{FFS}, the first term
\begin{displaymath}
  \sum_{l\geq 0}\frac{1}{l!}P_{l}(Y,\cdots,Y)
\end{displaymath}
in the above equation is equal to $(\hat{A}_{\hbar}\Ch)(Y)$. For the second
term, we can pull the sum into $\tr_{\gamma}$ and find
$\tr_{\gamma}(\exp_{\star}(Z))$. Hence,
\begin{displaymath}
  S(X)=(\hat{A}_{\hbar}\Ch)(Y)\Ch_{\gamma}(Z).
\end{displaymath}
Since $P_{k}^{\gamma}$ is the degree $k$ component of $S$, it is
equal to $(\hat{A}_{\hbar}\Ch\Ch_{\gamma})_{k}$.
\end{proof}
The same argument proves the twisted analogue for the cocycle
$\Theta^{V,\gamma}_{2k}$ defined in \eqref{Eq:theta}:
\begin{theorem}
\label{thm:local-trr}
  For $\dim V^\gamma \gg n$, the following identity holds:
  \begin{displaymath}
    [\ev_1\Theta_{2k}^{V,\gamma}]=(-1)^k
    \chi((\hat{A}_{\hbar}\Ch_V\Ch_{\gamma})_k).
  \end{displaymath}
\end{theorem}
\begin{remark}
In Theorem \ref{thm:local-trr}, $\Ch_V$ is a twisted Chern character
on $\mathfrak{gl}_V(\mathbb{K})^\gamma$ defined by
$\Ch_V(X):=\tr(\gamma \exp(X))$ for $X\in
\mathfrak{gl}_V(\mathbb{K})^\gamma$.
\end{remark}

%% file: algind.tex
\section{The algebraic index for orbifolds}
\label{Sec:algind}
In this section, we use the local Riemann-Roch formula in Theorem
\ref{thm:local-rr} to prove an algebraic index theorem on a
symplectic orbifold and thus confirm a conjecture by \cite{fst}.
As an application of this algebraic index theorem, we provide in
Section \ref{Sec:Kawasaki} an alternative proof of the
Kawasaki index theorem for elliptic operators on orbifolds \cite{Kaw:IEOVM}.
\subsection{The conjecture by Fedosov--Schulze--Tarkhanov}
Recall the set-up as described in Section \ref{Sec:preli}. Consider a
$G$-invariant formal deformation quantization $\calA^\hbar$ of $G_0$
with characteristic class $\Omega\in H^2(X,\K)$. Using the twisted trace
density of Section \ref{Sec:twhoch}, the trace
$\Tr:\calA^\hbar\rtimes G\rightarrow \K$ on the crossed product defined
by Equation \eqref{trace} gives rise to an index map by
\eqref{indexmap}. The following theorem computes this map in terms of
characteristic classes on the inertia orbifold $\tilde{X}$ of $X$.
\begin{theorem}
\label{thm:algind}
 Let $E$ and $F$ be orbifold vector bundles on $X$, which are isomorphic
 outside a compact subset. Then we have
\begin{displaymath}
\Tr_*([E]-[F])=
  \int_{\tilde{X}}\frac{1}{m}
  \frac{\Ch_\theta(\frac{R^E}{{2\pi i}}-\frac{R^F}{2\pi i})}{
  \det(1-\theta^{-1}\exp(-\frac{R^{\perp}}{2\pi i}))}
  \hat{A}\Big(\frac{R^T}{2\pi i}\Big)\exp\Big(
  -\frac{\iota^*\Omega}{2\pi i\hbar}\Big).
\end{displaymath}
\end{theorem}
\begin{proof}
First, observe that by definition the left hand side computes the pairing
between the cyclic cocycle of degree $0$ given by the trace, and $K$-theory.
Therefore, the left hand side only depends on the $K$-theory class of
$[E]-[F]\in K_{\text{\tiny orb}}^0(X)$ and will not change when we add a
trivial bundle to $E$ and $F$. Indeed, the right hand side is invariant
under such changes as well, as for a trivial bundle we have
$\Ch_\theta=\Ch=1$. Consequently, we can assume without loss
of generality that $\operatorname{rk} (E)=\operatorname{rk} (F)\gg n$,
enabling the use of the local Riemann--Roch Theorem \ref{thm:local-rr} of
the previous section.

Let us first compute the local index density of Theorem \ref{thm:local-rr} of
the trivial vector bundle. In this local computation we put $\theta=\gamma$ as the twisting automorphism. Notice that since $A$ is $G$-invariant,
$\iota^*A$ defined in Section 4.2, (\ref{pb}) is also $G$-invariant.
Therefore,
\begin{displaymath}
  \psi_{D}(1)=
  \ev_1(\Theta_{2k}^{N,\gamma})(\iota^*A\wedge\cdots\wedge \iota^*A)
\end{displaymath}
can be identified with the help of Theorem \ref{thm:local-rr} as follows.
Since the symplectic connection and $\iota^*A$ are
$\gamma$-invariant, it follows that the curvature $\iota^*\tilde{R}$
is also $\gamma$-invariant. As in \cite[Sect. 4.7]{FFS}, we thus get
\begin{displaymath}
  F(\iota^*A(\xi),\iota^*A(\eta))=\iota^*\tilde{R}(\xi,\eta)-
  \iota^*\Omega(\xi,\eta)
\end{displaymath}
for vector fields $\xi$ and $\eta$ on $B_{0}$. In the formula
above, $\tilde{R}\in\Omega^2(G_0,\mathfrak{sp}^\gamma_{2n})$ is
the curvature of $\nabla$, and we have used Equation
\eqref{pbcurv} to arrive at this result. To apply Theorem
\ref{thm:local-rr}, we split the curvature
\begin{displaymath}
  \iota^*\tilde{R}=\iota^*\tilde{R}^t+\iota^*\tilde{R}^{\perp}
\end{displaymath}
according to the decomposition
$\mathfrak{sp}^\gamma_{2n}=\mathfrak{sp}_{2k}\oplus
\mathfrak{sp}_{2(n-k)}^\gamma.$
Substituting this into the formula for
$\ev_1(\Theta_{2k}^{\gamma})$ of Theorem \ref{thm:local-rr}, we find
\begin{displaymath}
\begin{split}
\frac{(-1)^k}{(2k)!}\ev_1\Theta_{2k}^{\gamma}(\iota^*A ^{2k})& =
  (-1)^k \ev_1\Theta_{2k}^{\gamma}(\iota^*A\wedge\cdots\wedge \iota^*A)\\
  &=\frac{1}{n!}P^{\gamma}_k\big((\iota^*\tilde{R}-\iota^*\Omega)^k\big)\\
  &=\big(\hat{A}_{\hbar}(\iota^*\tilde{R}^T)\Ch(-\iota^*\Omega)
  \Ch_{\gamma}(\iota^* \tilde{R}^\perp)\big)_k\\
  &=\hbar^k
  \Big( \hat{A}(R^T)\Ch\Big(-\frac{1}{\hbar}\iota^*\Omega\Big)\tr_{\gamma}
  \Big(\exp_{\star}\Big(\frac{\iota^*\tilde{R}^\perp}{\hbar}
  \Big)\Big)\Big)_{k}.
\end{split}
\end{displaymath}

By the constructions in \cite[Sec.~5]{fe:g-index} and the computations in
\cite[Thm.~4.1]{fe:g-index}, one has
$\tr_{\gamma}(\exp_{\star}(\frac{\iota^*\tilde{R}^\perp}{\hbar}))=
(\det(1-\gamma^{-1}\exp(-R^{\perp})))^{-1}$.
Inserting this expression into the above formula, we obtain
\begin{displaymath}
  \psi_D(1)=\ev_1(\Theta_{2k}^{\gamma})=
  \hbar^k\frac{\hat{A}(R^T)\Ch(-\frac{\iota^*\Omega}{\hbar})}{
  \det(1-\gamma^{-1}\exp(-R^{\perp}))}.
\end{displaymath}
For a general (i.e., nontrivial) vector bundle $E$, we put $V$ equal to the
typical fiber of $E$ and use Theorem \ref{thm:local-trr} to find
\begin{displaymath}
  \psi_{D_E}(1)=\ev_{1}(\Theta_{2k}^{V,\gamma})=
  \hbar^{k}\frac{\Ch_{\gamma}(R^E)}{\det(1-\gamma^{-1}\exp(-R^{\perp}))}
  \hat{A}\Big(R^T\Big)\Ch\Big(-\frac{\iota^*\Omega}{\hbar}\Big).
\end{displaymath}
Applying Proposition \ref{prop:k}, the result now follows.
\end{proof}
\subsection{The Kawasaki index theorem}
\label{Sec:Kawasaki}
In this section, we derive Kawasaki's index theorem
\cite{Kaw:IEOVM} from Theorem \ref{thm:algind} for orbifolds.
Hereby, we apply the methods introduced by Nest--Tsygan in
\cite{nt96}, where the relation between formal and analytic index
formulas has been studied for compact riemannian manifolds.

Let $X$ be a reduced compact riemannian orbifold, represented by a
proper \'etale Lie groupoid $G_1\rightrightarrows G_0$. We look at
its cotangent bundle $T^* X$, which is  a symplectic orbifold
and is represented by the groupoid $T^*G_1\rightrightarrows T^*G_0$ with
the canonical (invariant) symplectic structure $\omega$ on $T^*G_0$.
The standard asymptotic calculus of pseudodifferential operators on
manifolds can be extended to $X$, as explained in Appendix \ref{Sec:appb}.

According to (\ref{Eq:AsExOpStar}), the operator
product on pseudodifferential operators defines an invariant star
product $\star^\text{\tiny\rm Op}$ on $T^*G_0$. Since it is invariant,
this star product descends to a star product $\star^\text{\tiny\rm Op}_X$
on $T^*X$. By ({\ref{Eq:OpTr}}) and (\ref{Eq:InvOpTr}), the operator trace
on trace class operators defines a trace $\Tr^\text{\tiny\rm Op}$ on the
deformation quantization of $T^*G_0$ and also a trace
$\Tr^\text{\tiny\rm Op}_X$ on $T^*X$.

In the case of a reduced orbifold $T^*X$, one knows that the
deformed groupoid algebra on $T^*G$ constructed by Fedosov's
approach as in \cite{ta:quantization} is Morita equivalent to the
deformed $G$-invariant algebra on $T^*G_0$ with $\star^\text{\tiny\rm Fe}$
as product \cite[Prop.6.5]{nppt}.
Recall that the latter is a deformation quantization for $T^*X$.
By this Morita equivalence, $\Tr$ as defined
by (\ref{trace}) reduces to a trace $\Tr^\text{\tiny\rm Fe}_X$ on the
deformation quantization $\star^\text{\tiny\rm Fe}_X$ of $T^*X$.

By the classification of star products on a symplectic manifold by
Nest and Tsygan \cite{nt}, the deformation quantization on
$T^*G_0$ constructed in Appendix \ref{Sec:appb} is equivalent to
the $G$-invariant deformation quantization $\star^\text{\tiny\rm
Fe}$ of $T^*G_0$ obtained by Fedosov quantization with Weyl
curvature equal to $-\omega$. Since $G$ is a proper \'etale Lie
groupoid, one can choose the equivalence  $\Phi$ between
$\star^\text{\tiny\rm Op}$ and $\star^\text{\tiny\rm Fe}$ to be
$G$-invariant. The pull-back of the trace $\Tr^\text{\tiny\rm
Fe}_X$ by  $\Phi$ thus defines  a trace on $\star^\text{\tiny\rm
Op}_X$. We denote  this pull-back trace by $\Tr^\text{\tiny\rm
Fe}_X$ as well.
\begin{proposition}
\label{prop:optr=fetr}
  The pull-back trace $\Tr^\text{\tiny\rm Fe}_X$ is equal
  to the trace $\Tr^\text{\tiny\rm Op}_X$.
\end{proposition}
\begin{proof}
Consider a smooth function $f$ on $T^*X$ with support in an open subset
$O \subset T^*X$ such that $O$ is isomorphic to the orbit space of a
transformation groupoid $\tilde{O}\rtimes \Gamma$, where
$\tilde{O}$ is a $\Gamma$-invariant open subset of some
finite dimensional symplectic vector space on which $\Gamma$ acts by
linear symplectomorphisms.
When restricted to $O$, the traces $\Tr^\text{\tiny\rm Fe}_O$ and
$\Tr^\text{\tiny\rm Op}_O$ both correspond to operator traces, but
are defined by different polarizations of $\tilde{O}$. By construction,
$\Tr^\text{\tiny\rm Op}_X$ is defined by a real polarization, while
$\Tr^\text{\tiny\rm Fe}_X$ is defined by a complex polarization.
The Hilbert spaces corresponding to these different polarizations are
related by a $\Gamma$-invariant unitary operator. Hence, one concludes that
$\Tr^\text{\tiny\rm Fe}_O=\Tr^\text{\tiny\rm Op}_O$, and therefore
$\Tr^\text{\tiny\rm Fe}_X(f)=\Tr^\text{\tiny\rm Op}_X(f)$.

To check that $\Tr^\text{\tiny\rm Fe}_X(f)=\Tr^\text{\tiny\rm Op}_X(f)$
for any compactly supported smooth function $f$ on $T^*X$ one now uses
an appropriate partition of unity to reduce the claim to the local case
which just has been proved. The proposition follows.
\end{proof}

The symbol $\sigma(D)$ of an elliptic operator $D$ on $X$ defines
an isomorphism between two orbifold vector bundles $E$ and $F$ on
$T^*X$ outside a compact subset. Similar to
\cite[Eq.~(4.2.2)]{fe:book}, the index of $D$ can be computed by
$\Tr^\text{\tiny\rm Op}_X([E]-[F])$. By Proposition \ref{prop:optr=fetr},
we can use $\Tr^\text{\tiny\rm Fe}_X([E]-[F])$ to calculate
this. By  Theorem \ref{thm:algind} we have
\begin{displaymath}
 \Tr^\text{\tiny\rm Fe}_X([E]-[F]) =
 \int_{\widetilde{T^*X}}\frac{1}{m}\frac{\Ch_\theta
 (\frac{R^E}{2\pi i}-\frac{R^F}{2\pi i})}{\det(1-\theta^{-1}
 \exp(-\frac{R^{\perp}}{2\pi i}))} \hat{A}\Big(\frac{R^T}{2\pi i}
 \Big)\Ch\Big(-\frac{\iota^*\Omega}{2\pi i\hbar}\Big).
\end{displaymath}
Now recall that in the construction of the star product
$\star^\text{\tiny\rm FE}_X$, we have fixed $-\iota^*\Omega$ to
$\iota^*\omega$. Since the canonical symplectic form $\omega$ is
exact, one concludes that $\Ch(-\frac{1}{\hbar}\iota^*\Omega)$ is
equal to $1$ in this case. Hence we obtain the following theorem
as an application of Theorem \ref{thm:algind}.
\begin{theorem}\cite{Kaw:IEOVM}
\label{thm:kawasaki}
  Given an elliptic operator $D$ on a reduced compact orbifold $X$, one has
\begin{displaymath}
 \operatorname{index} (D)=\int_{\widetilde{T^*X}}\frac{1}{m}
 \frac{\Ch_\theta\big(\frac{\sigma(D)}{2\pi i}\big)}{
 \det \big(1-\theta^{-1}\exp \big(-\frac{R^{\perp}}{2\pi i}\big)\big)}
 \hat{A}\Big(\frac{R^T}{2\pi i}\Big),
\end{displaymath}
where $\sigma(D)$ is the symbol of $D$.
\end{theorem}

%% file: appa.tex
\section{Twisted Hochschild and Lie algebra cohomology}
\label{Sec:appa}
In this section we consider the cohomology of the Lie algebra
$\gl_N (  \W _{2n } \rtimes \Gamma)$ for $N\gg 0$, where $ \W _{2n
}$ is the (formal) Weyl algebra over $\R^{2n}$ and $\Gamma$ is a
finite group which is assumed to act effectively by
symplectomorphisms on $\R^{2n}$. More precisely, we will compute
for $0\leq p \leq n_\Gamma$, $N\gg n$, and $q\in \N^*$ the Lie
algebra cohomology groups
\begin{equation}
\label{EqLieAlgHomGps}
   H^p \big( \gl_N ( \W _{2n }\rtimes \Gamma ) ;
   S^q \mathfrak{M}_N ( \W _{2n } \rtimes \Gamma)^* \big),
\end{equation}
where $n_\Gamma$ is a natural number depending on $\Gamma$ and
$S^q M$ denotes for every vector space $M$ the $q$-th symmetric
power. Clearly, if $M$ is a bimodule over an algebra $A$, then
$S^qM$ carries the structure of a $\gl_N (A)$-module as follows:
\begin{displaymath}
  [a, m_1 \vee \ldots \vee m_q] = \sum_{i=1}^q m_
  1 \vee \ldots \vee
  [a,m_i]\vee \ldots \vee m_q ,
\end{displaymath}
where $a \in A$, $m_1, \ldots ,m_q\in M$ and $[a,m_i] = a \, m_i -m_i \, a$.
\subsection{$\Z_2$-graded Hochschild and Lie algebra cohomology}
For the computation of the above Lie algebra cohomology we will make use
of the super or in other words $\Z_2$-graded versions of Hochschild and
Lie algebra (co)homology. Let us briefly describe the construction of these
super homology theories. To this end consider first a $\Z_2$-graded unital
algebra $A=A_0\oplus A_1$. For a homogeneous element $a\in A$ we then
denote by $|a|$ its degree, that means the unique element $i\in \Z_2$
such that $a \in A_i$.
Now there exist uniquely defined face maps
$b_j : C_p(A) \rightarrow C_{p-1}(A)$, $0\leq j \leq p$ which satisfy
the following relations for homogeneous
$a_0, a_1 ,\ldots, a_p\in A$:
\begin{equation}
\begin{split}
  b_j (a_0 \otimes & \ldots \otimes a_p) = \\
  & = \begin{cases}
     a_0 \otimes \ldots \otimes a_i a_{i+1} \otimes \ldots \otimes a_p, &
     \text{for $0\leq j < p$},\\
     (-1)^{|a_p| \, (|a_0| + \ldots + |a_{p-1}|)}
     a_p a_0 \otimes a_1 \otimes \ldots \otimes a_{p-1}, &
     \text{for $j = p$}.
  \end{cases}
\end{split}
\end{equation}
Moreover, one has degeneracy maps $s_j : C_p(A) \rightarrow C_{p+1}(A) $,
$0\leq j\leq p$,
and cyclic operators $t_p : C_p(A) \rightarrow C_p(A) $ defined as
follows:
\begin{equation}
\begin{split}
  s_j & (a_0 \otimes  \ldots \otimes a_p) = a_0 \otimes \ldots \otimes a_j
  \otimes 1 \otimes a_{j+1} \otimes \ldots \otimes a_p , \\
  t_p & (a_0 \otimes  \ldots \otimes a_p) =
  (-1)^{|a_p| \, (|a_0| + \ldots + |a_{p-1}|)}
  (a_p \otimes a_0 \otimes  \ldots \otimes a_{p-1}) .
\end{split}
\end{equation}
It is straightforward to check that these data give rise to a cyclic object
in the category of super vector spaces, hence to the Hochschild
and cyclic homology of the super algebra $A$.

As a particular example consider the super algebra $\C [\varepsilon]$, where
$\epsilon$ is of degree $1$ and satisfies $\varepsilon^2 =0$
(cf.~\cite[Sec.~3.1]{FeiTsy:RRTLACI}). One
proves immediately, that for each $p$, the chains
$\omega^1_p := 1\otimes \varepsilon \otimes \ldots \otimes \varepsilon $ and
$\omega^2_p := \varepsilon \otimes \varepsilon \otimes \ldots \otimes
\varepsilon $
are Hochschild cycles which generate the $\Z_2$-graded Hochschild homology
$HH_p (\C[\varepsilon])$. Moreover, one checks that
$B ( \omega^1_p) = B (\omega^0_p) =0$ and $B(\omega^2_p )=\omega^1_{p+1}$,
where $\omega^0_p$ denotes the cycle $1\otimes \ldots \otimes 1$.
Consider Connes' periodicity exact sequence
\begin{displaymath}
  \longrightarrow HH_p (\C[\varepsilon] ) \overset{I}{\longrightarrow}
  HC_p (\C[\varepsilon] ) \overset{S}{\longrightarrow}
  HC_{p-2} (\C[\varepsilon] ) \overset{B}{\longrightarrow}
  HH_{p-1} (\C[\varepsilon] ) \longrightarrow ,
\end{displaymath}
which also holds in the $\Z_2$-graded version.
Then one checks by induction on $p$ that $HC_p(\C[\varepsilon] )$
is generated by $[\omega^0_p]$ and $[\omega^2_p]$ in case $p$ is
even, and by $[\omega^0_p]$, if $p$ is odd. Thus, one has in particular
\begin{equation}
  HC_p (\C[\varepsilon] ) \cong HC_p (\C) \oplus \C \quad
  \text{for every $p\in \N$}.
\end{equation}

Next consider a super Lie algebra $\mathfrak g = \mathfrak g_0
\oplus \mathfrak g_1$. The super Lie algebra homology or in other
words $\Z_2$-graded Lie algebra homology $H_\bullet (\mathfrak g ,
\C)$ is then defined as the homology of the complex
\begin{displaymath}
 \longrightarrow \Lambda^l \mathfrak g \overset{d}{\longrightarrow}
 \Lambda^{l-1} \mathfrak g \overset{d}{\longrightarrow} \ldots
 \overset{d}{\longrightarrow}
 \Lambda^1 \mathfrak g \overset{d}{\longrightarrow}
 \Lambda^0 \mathfrak g =\C ,
\end{displaymath}
where $\Lambda^l \mathfrak g := \bigoplus_{p+q=l} E^p \mathfrak g_0
 \oplus S^q \mathfrak g_1$ means the super exterior product, and
\begin{displaymath}
\begin{split}
 d (\xi_1\wedge \ldots \wedge \xi_p) =
 \sum_{1\leq i < j\leq p}
 [\xi_i ,\xi_j] \wedge \xi_1 \wedge \ldots \wedge \hat \xi_i \wedge
 \ldots \wedge \hat \xi_j \wedge \ldots \wedge \xi_p.
\end{split}
\end{displaymath}
Clearly, every super algebra $A$ gives rise to a super Lie algebra
$\mathfrak {gl}_N (A)$. Recall that one has the following relation
between the (super) Lie algebra homology and the $\Z_2$-graded
cyclic homology of $A$ (cf.~\cite[Thm.~10.2.5]{LoCH}):
\begin{equation}
\label{lodquil}
  H_\bullet (\mathfrak{gl}_N (A) ;\C) \cong \Lambda_\bullet
  (HC[1](A)),\ \ \ \ N\gg 0,
\end{equation}
where $\Lambda_\bullet$ denotes the functor which associates to a graded
vector space its graded symmetric algebra.

\subsection{Computation of the Lie algebra cohomology}
Let us come back to our original goal, the computation of the
Lie algebra homology groups in (\ref{EqLieAlgHomGps}).
To this end recall first the following result which has been proved in
various forms in \cite{AFLS:HIAWAGF,nppt,dolget}.
\begin{proposition}
\label{Prop:HHConv}
  Let $\Gamma$ be a finite group which acts (from the right) by
  automorphisms on an algebra $A$ over a field $\Bbbk$. Then the Hochschild
  (co)homology of the convolution algebra $A\rtimes \Gamma$ statisfies
  \begin{align}
    HH_\bullet (A\rtimes \Gamma) & = H_\bullet (A, A\rtimes \Gamma)^\Gamma, \\
    HH^\bullet (A\rtimes \Gamma) & = H^\bullet (A, A\rtimes \Gamma)^\Gamma
  \end{align}
\end{proposition}
Now let us consider again the (formal) Weyl algebra $ \W _{2n }$ on
$\R^{2n}$ (over the field $\K=\C((\hbar))$) and let $\Gamma$ act by
symplectomorphisms on $\R^{2n}$. Then note that $H_p \big(\gl_N ( \W
_{2n }\rtimes \Gamma); S^q \mathfrak{M}_N ( \W _{2n } \rtimes
\Gamma)\big)$ is dual to $H^p (\gl_N ( \W _{2n }\rtimes \Gamma );
S^q \mathfrak{M}_N ( \W _{2n } \rtimes \Gamma)^*)$, the Lie algebra
cohomology we are interested in. Hence it suffices to determine the
homology groups $H_p \big(\gl_N ( \W _{2n }\rtimes \Gamma); S^q
\mathfrak{M}_N ( \W _{2n } \rtimes \Gamma)\big)$. Define the super
algebra $A$ as the $\Z_2$-graded tensor product
\begin{displaymath}
  A = ( \W _{2n } \rtimes \Gamma) \otimes_\C \C[\varepsilon] .
\end{displaymath}
By the above definition of super Lie algebra homology it is now clear that
\begin{displaymath}
  H_p ( A,\K) \cong \bigoplus_{p+q=k} H_p ( \W _{2n } \rtimes
  \Gamma;
  S^q  \W _{2n } \rtimes \Gamma).
\end{displaymath}
Recall from Prop.~\ref{computationtwh} that for every $\gamma \in
\Gamma$ the (twisted) Hochschild  homology $H^\gamma_\bullet ( \W
_{2n }) =H_\bullet ( \W _{2n }, {\W _{2n }}_\gamma)$ is
concentrated in degree $2k_{\langle\gamma\rangle}$, where $2
k_{\langle\gamma\rangle}$ is the dimension of the fixed point
space of $\gamma$ (which depends only on the conjugacy class
$\langle\gamma\rangle$). Using the twisted Hochschild (co)homology
of $ \W _{2n }$ and Prop.~\ref{Prop:HHConv} one can derive the
following result exactly as in \cite{AFLS:HIAWAGF}, where the case
of the nonformal Weyl algebra has been considered.
\begin{theorem}
\label{Thm:DimHCon}
  {\rm (}cf.~\cite[Thm.~6.1]{AFLS:HIAWAGF}{\rm )}
  Let $\Gamma$ be a finite group which acts effectively by
  symplectomorphisms on $\R^{2n}$ and consider its induced action (from the
  right) on the formal Weyl algebra over $\R^{2n}$. Denote for every
  $p \in \N$ by $l_p (\Gamma)$ the number of conjugacy classes of elements
  of $\Gamma$ having a $p$-dimensional fixed point space. Then the following
  formula holds for the Hochschild (co)homology of the crossed product algebra
  $ \W _{2n } \rtimes \Gamma$:
  \begin{equation}
     \dim_\K H_p ( \W _{2n } \rtimes \Gamma) =
     \dim_\K H^{2n- p} (  \W _{2n } \rtimes \Gamma ,
     \W _{2n } \rtimes \Gamma ) = l_p (\Gamma).
  \end{equation}
\end{theorem}
Now choose for every element $\gamma$ of $\Gamma$ a Hochschild
cycle $\alpha_\gamma$ the homology class of which generates $H_{2
k_{\langle\gamma\rangle}}( \W _{2n }, {\W _{2n }}_\gamma)$, and
denote by $\sqcup$ the exterior (shuffle) product. Then it is
clear that
\begin{equation}
\label{Eq:quasiisomixed}
  (C_{\bullet - 2k_{\langle\gamma\rangle}} (\K),b,B)
  \underset{\sqcup \alpha_\gamma}{\longrightarrow}
  (C_\bullet ( \W _{2n }),b_{\gamma} ,B)
\end{equation}
is a morphism of mixed complexes. Moreover, by the properties of $
\W _{2n }$ it is even a quasi-isomorphism for twisted Hochschild
homology, hence a quasi-isomorphism of mixed complexes. Next
consider the following composition of morphisms of mixed complexes
\begin{equation}
\begin{split}
  \bigoplus_{\langle\gamma\rangle \in \operatorname{Conj}(\Gamma)}
  \big( C_{\bullet -2k_{\langle\gamma\rangle}} (\K),b,B \big) &
  \overset{\alpha}{\rightarrow}
  \big( C_\bullet ( \W _{2n }, \W _{2n }\rtimes \Gamma),b,B \big)
  \overset{\iota}{\rightarrow}
  \big( C_\bullet ( \W _{2n }\rtimes \Gamma),b,B \big)
  \\
  (c_{\langle\gamma\rangle})_{\langle\gamma\rangle\in\operatorname{Conj}
  (\Gamma)}
  &\mapsto (c_{\langle\gamma\rangle} \sqcup \alpha_\gamma
  )_{\gamma \in \Gamma} ,
\end{split}
\end{equation}
where $\operatorname{Conj}(\Gamma)$ denotes the set of conjugacy classes of
$\Gamma$, and where we have used the natural identification
\begin{displaymath}
   \big( C_\bullet ( \W _{2n }, \W _{2n }\rtimes \Gamma),b,B \big) \cong
  \bigoplus_{\gamma \in \Gamma} \big( C_\bullet ( \W _{2n }),
  b_\gamma, B \big) .
\end{displaymath}
One now checks easily that the $\alpha_\gamma$ can be chosen in such a way that
\begin{equation}
  \alpha_{\gamma\tilde\gamma\gamma^{-1}} = \gamma\, \alpha_{\tilde\gamma}
 \quad \text{for all $\gamma, \tilde \gamma\in \Gamma$}.
\end{equation}
Under this assumption on the $\alpha_\gamma$ the
composition $\iota \circ \alpha$ is a quasi-isomorphism of mixed complexes.
Let us show this in some more detail.
By the choice of the $\alpha_\gamma$ it is clear that the image of the
morphism $\alpha$ is invariant under the action of $\Gamma$. Moreover,
$\alpha$ is injective on homology by construction.
Since the restriction of $\iota$ to the invariant part is a quasi-isomorphism
in Hochschild homology by  Prop.~\ref{Prop:HHConv}, a simple dimension
counting argument then shows by Thm.~\ref{Thm:DimHCon} that the composition
$\iota \circ$ is a  quasi-isomorphism in Hochschild homology as well.
Hence $\iota \circ \alpha$ is also a quasi-isomorphism of mixed complexes.
Thus, if $2k$ is minimal among the dimensions of fixed point spaces of the
elements of $\Gamma$, one gets
\begin{equation}
  HC_p ( \W _{2n } \rtimes \Gamma ) =
  \begin{cases}
    0 & \text{for $0 \leq p < 2 k$}, \\
    \K^{l_{2k}(\Gamma)} & \text{for $p = 2 k$}.
  \end{cases}
\end{equation}

By construction and the Kunneth-isomorphism, the Hochschild complex
$C_\bullet (A)$ is quasi-isomorphic to the tensor product complex
$(C_\bullet (A \rtimes \Gamma ),b) \otimes C_\bullet (\C [\varepsilon])$.
Hence, by the above considerations
\begin{displaymath}
  \bigoplus_{\langle\gamma\rangle \in \operatorname{Conj} (\Gamma)}
  \big( C_{\bullet -k_{\langle\gamma\rangle}} (\K [\varepsilon]),b,B\big)
  \longrightarrow
  (C_\bullet (A),b ,B)
\end{displaymath}
is a quasi-isomorphism of mixed complexes as well, which implies that
\begin{equation}
  HC_p (A) =
  \begin{cases}
    0 & \text{for $0 \leq p < 2k$}, \\
    \K^{2l_{2k}(\Gamma)} & \text{for $p = 2k $}.
  \end{cases}
\end{equation}
Inspecting formula \eqref{lodquil} both for $ \W _{2n }$ and $A$
one now obtains the main result of this section.
\begin{theorem}
\label{Thm:LieAlgHomConv} Let $\gamma$ be a linear
symplectomorphism of $\R^{2n}$ of finite order, $\Gamma$ be the
cyclic group generated by $\gamma$, and put $2k := \dim (\R^{2n})^\gamma$.
Then for $N\gg n$,
\begin{equation}
  H_p (\gl_N( \W _{2n } \rtimes \Gamma) ;S^q \mathfrak{M}_N
  ( \W _{2n } \rtimes \Gamma)) =
  \begin{cases}
    0  & \text{for $0 \leq p  < 2 k$}, \\
    \K^l & \text{for $p = 2k $},
  \end{cases}
\end{equation}
where $l:= l_{2k} (\Gamma)$ is the number of elements of $\Gamma$
having a fixed point space of dimension $2k$.
\end{theorem}
\begin{remark}
  In the ``untwisted case'', i.e.~for $\Gamma$ the trivial group, these
  Lie algebra homology groups have been determined in
  \cite{FeiTsy:RRTLACI}.
\end{remark}
  By Eq.~(\ref{lodquil}) one concludes easily that the Lie algebra
  (co)homology groups of $\gl_N ( \W _{2n } \rtimes \Gamma)$ with values
  in $S^q \mathfrak{M}_N ( \W _{2n } \rtimes \Gamma)$ are
  Morita invariant. But one knows that $ \W _{2n }\rtimes \Gamma$ is Morita
  equivalent to the invariant algebra $ \W _{2n }^\Gamma$ (see \cite{dolget}).
  Thus the preceding theorem entails immediately
\begin{corollary}
\label{Cor:MainApp}
 Under the assumptions from above, one has for the invariant algebra
 $ \W _{2n }^\gamma$,
 \begin{equation}
  H_p (\gl_N ( \W _{2n }^\gamma) ;S^q \mathfrak{M}_N ( \W _{2n }^\gamma)) =
  \begin{cases}
    0  & \text{for $0 \leq p  < 2 k$}, \\
    \K^l & \text{for $p = 2k $}.
  \end{cases}
 \end{equation}
\end{corollary}

We can extend the above Corollary \ref{Cor:MainApp} to the Weyl
algebra with twisted coefficients
$\W_{2n}^V:=\W_{2n}\otimes \End(V)$,
where $V$ is a complex vector space $\gamma$ acts on diagonally.

Note that $\W_{2n} ^V$ is Morita equivalent to $\W_{2n}$ by
the bimodule $\W_{2n}\otimes V$, and $\W_{2n} ^V\rtimes \Gamma$ is
Morita equivalent to $\W_{2n}\rtimes \Gamma$ by the bimodule
$(\W_{2n}\otimes V)\rtimes \Gamma$. We have that
\[
\begin{array}{ll}
HH^{\bullet}(\W_{2n}^V\rtimes
\Gamma)&=HH^{\bullet}(\W_{2n}\rtimes \Gamma),\\
HH_{\bullet}(\W_{2n}^V\rtimes \Gamma)&=HH_{\bullet}(\W_{2n}\rtimes
\Gamma).
\end{array}
\]

Hence, by the same arguments as Corollary \ref{Cor:MainApp}, we
have the following result for the Lie algebra homology of
$\mathfrak{gl}_N((\W_{2n}^V )^{\gamma})$.
\begin{corollary}
\label{Cor:mainapp-twist} For $(\W_{2n}^V)^{\gamma}$, one has when
$N\gg n$,
 \begin{equation}
  H_p (\gl_N ((\W_{2n}^V)^\gamma) ;S^q \mathfrak{M}_N ((\W_{2n}^{V})^\gamma)) =
  \begin{cases}
    0  & \text{for $0 \leq p  < 2 k$}, \\
    \K^l & \text{for $p = 2k $}.
  \end{cases}
 \end{equation}
\end{corollary}

%% file: appb.tex
\section{Asymptotic pseudodifferential calculus}
\label{Sec:appb}
In the following we sketch  the construction of the asymptotic calculus
of pseudodifferential operators on a proper \'etale Lie groupoid. Hereby,
we adapt the presentation in \cite{nt96} to the groupoid case. For more
details on the pseudodifferential calculus for proper \'etale groupoids, we
refer the reader to  Hu's thesis \cite{hu:thesis}. See also
\cite{Wid:CSCPO} for the complete symbol calculus on riemannian manifolds,
and \cite{P1} for its application to deformation quantization.

Let $G \rightrightarrows G_0$ be a proper \'etale Lie groupoid
and fix a $G$-invariant riemannian metric on $G_0$.
Then check that the
canonical action of $G$ on the diagonal
$\Delta \subset G_0\times G_0$ can be extended to a $G$-action on a
neighborhood $W\subset G_0\times G_0$ of $\Delta$, since $G$ is proper
\'etale. Next construct a cut-off function
$\chi :G_0\times G_0 \rightarrow [0,1]$ with the following properties:
\begin{enumerate}
\item
  $\supp \chi \subset W$, and the restriction $\chi_{|W}$ is invariant
  under the diagonal action of $G$ on $W$.
\item
  One has $\chi(x,y)=\chi(y, x)$ for all $x,y\in G_0$.
\item
  On a neighborhood of the diagonal, one has $\chi \equiv 1$.
\item
  For every $x\in G_0$, the set
  $Q_x=\{ y \in G_0 \mid (x,y)\in \supp (\chi) \}$ is compact and
  geodesically convex with respect to the chosen metric.
\end{enumerate}

For every open $U\subset G_0$ denote by $\Sym^m(U)$, $m\in \Z$, the
space of symbols of order $m$ on $U$, that means the space of
smooth functions $a$ on $T^*U$ such that in each local coordinate system
of $U$ and each compact set $K$ in the domain of the local coordinate
system there is an estimate of the form
\begin{displaymath}
  \big| \partial_x^{\alpha}\partial_{\xi}^{\beta} a (x,\xi)\big| \leq
  C_{\alpha, \beta}(1+|\xi|^2)^{\frac{m-|\beta|}{2}}, \quad
  \text{$x\in K$, $\xi \in T^*_x G_0$, $\alpha,\beta \in \N^n$},
\end{displaymath}
for some $C_{K,\alpha,\beta} >0$. Clearly, the $\Sym^m(U)$ are the
section spaces of a sheaf $\Sym^m$ on $G_0$. Moreover, each
$\Sym^m$ is even a $G$-sheaf, i.e.~carries a right $G$-action on
the stalks. Finally, one obtains two further symbol sheaves by
putting
\begin{displaymath}
 \Sym^\infty := \bigcup_{m\in \Z}\Sym^m, \quad
 \Sym^{-\infty} := \bigcap_{m\in \Z}\Sym^m.
\end{displaymath}
Similarly, one constructs the presheaves $\PDO^m$
of pseudodifferential operators of order
$m\in \Z \cup \{ - \infty , \infty\}$ on $G_0$.
Next let us recall the definition of the symbol map
$\sigma $ and its quasi-inverse, the quantization map $\Op$.
The symbol map associates to every operator $A\in \PDO^m (U)$
a symbol $a \in \Sym^m (U)$ by setting
\begin{displaymath}
  a (x,\xi) : = A \big( \chi (\cdot, x)
  e^{i \langle \xi , \Exp_x^{-1} (\cdot )\rangle } \big) \, (x) ,
\end{displaymath}
where $\Exp_x^{-1}$ is the inverse map of the exponential map on $Q_x$.
The quantization map is given by
\begin{displaymath}
\begin{split}
  \Op: \: & \Sym^m (U) \rightarrow \PDO^m (U) \subset \Hom \big(
  \mathcal{C}^{\infty}_\text{\tiny\rm cpt} (U),
  \mathcal{C}^{\infty}(U)\big), \\
  & \big( \Op (a) f\big) (x) :=
  \int_{T^*_x G_0} \int_{G_0} e^{i \langle \xi , \Exp_x^{-1} (y)\rangle }
  \chi(x, y)  a (x , \xi) f(y) \, dy\,  d\xi .
\end{split}
\end{displaymath}
The maps $\sigma$ and $\Op$ are now quasi-inverse to each other in the sense
that the induced morphisms $\overline{\sigma}$ and $\overline{\Op}$
between the quotient sheaves $\Sym^\infty / \Sym^{-\infty}$
and $\PDO^\infty / \PDO^{-\infty}$ are isomorphisms such that
$\overline{\Op}^{-1} = \overline{\sigma} $.

By the space $\ASym^m (U)$, $m\in \Z$ of {\it asymptotic symbols} over an open
$U\subset G_0$ one understands the space of all
$q \in \calC^\infty (T^*U \times [0 ,\infty))$ such that for each
$\hbar \in [0 , \infty)$ the function $q(-,\hbar)$ is in $\Sym^m (U)$
and such that $q$ has an asymptotic expansion of the form
\begin{displaymath}
  q \sim \sum_{k\in \N} \hbar^k a_{m-k} ,
\end{displaymath}
where each $a_{m-k}$ is a symbol in $\Sym^{m-k} (U)$. More precisely,
this means that one has for all $N \in \N$
\begin{displaymath}
  \lim_{\hbar \searrow \, 0} \Big( q (-,\hbar)  - \hbar^{-N}
  \sum_{k=0}^N \hbar^k a_{m-k}\Big) = 0 \quad \text{in $\Sym^{m-N} (U)$}.
\end{displaymath}
Clearly, the $\ASym^m (U)$ are the sectional spaces of a $G$-sheaf $\ASym^m$.
By forming the union resp.~intersection of the sheaves $\ASym^m$ like above
one obtains two further $G$-sheaves $\ASym^\infty$ and $\ASym^{- \infty}$.
Moreover,  since $G$ acts on these in a natural way, we also obtain the sheaf
$\ASym^{\pm \infty}_G := (\ASym^{\pm \infty})^G$ of invariant asymptotic
symbols and the convolution algebra $\ASym^\infty \rtimes G$.
For $m\in \Z \cup \{ -\infty,\infty\}$ consider now the subsheaves
$\JSym^m  \subset \ASym^m$ and $\JSym^m_G \subset \ASym^m_G$ consisting of
all (invariant) asymptotic symbols which vanish to infinite order on
$\hbar =0$. The quotient sheaves  $\A^m := \ASym^m / \JSym^m$ and
$\A_G^m := \ASym_G^G / \JSym_G^m$
can then be identified with the formal power
series sheaves $\Sym^m [[\hbar]]$ resp.~$\Sym^m_G [[\hbar]]$.

The operator product on $\PDO^\infty$ now induces a
(asymptotically associative) product on $\ASym^\infty (G_0)$ by
defining for $q,p \in \ASym^\infty (G_0)$
\begin{equation}
   q \star p  (-,\hbar) :=
   \begin{cases}
      \iota_{\hbar^{-1}} \sigma \big( \Op ( \iota_\hbar q (-,\hbar) ) \circ
      \Op (\iota_\hbar p (-,\hbar )) \big) & \text{if $\hbar > 0$} ,\\
      q(-,\hbar ) \cdot p (-,\hbar )& \text{if $\hbar =0$}.
   \end{cases}
\end{equation}
Hereby, $\iota_\hbar : \Sym^\infty (G_0) \rightarrow \Sym^\infty (G_0)$
is the map which maps a symbol $a$ to the symbol
$(x,\xi) \mapsto a(x,\hbar \xi)$.
By standard techniques of pseudodifferential calculus (cf.~\cite{P1}),
one checks that the ``star product'' $\star$ has an asymptotic expansion
of the following form:
\begin{equation}
\label{Eq:AsExOpStar}
  q \star p \sim q \cdot p + \sum_{k =1}^\infty c_k (q,p) \, \hbar^k ,
\end{equation}
where the $c_k$ are bidifferential operators on $T^*G_0$ such that
\begin{displaymath}
    c_1 (a,b) -c_1 (b,a) = -i \{a,b\} \quad
    \text{for all symbols $a,b \in \Sym^\infty (G_0)$}.
\end{displaymath}
Hence, $\star$ is a star product on the quotient sheaf $\A^\infty$
and also on $\A_G^\infty$, since by construction the star
product of invariant symbols is again invariant. Thus one obtains
deformation quantizations for both the sheaf $\calA_{T^*G_0}$ of smooth
functions on $TG_0$ and the sheaf $\calA_{T^*X} = \calA^G_{T^*G_0}$
of smooth functions on the orbifold $X$ represented by the groupoid $G$.

The invariant riemannian metric on $G_0$ gives rise to Hilbert spaces
$L^2 (G_0)$ and $L^2 (X) = \pi_G L^2 (G_0)$, where $\pi_G$ is the
orthogonal projection on the space of invariant functions.
Hence there is a natural operator trace $\Tr_{L^2}$
on the space $\PDO^{- \dim X}_\text{\tiny \rm cpt} (G_0)$ of
pseudodifferential operators of order $\leq -\dim X$ with compact support.
Thus there is a map
\begin{displaymath}
  \Tr^\text{\tiny \rm Op} : \A^{- \infty}_\text{\tiny \rm cpt} (G_0)
  \rightarrow \C[\hbar^{-1} ,\hbar]], \quad
  q \mapsto \Tr_{L^2} \big( \Op \iota_\hbar (q(-,\hbar)) \big),
\end{displaymath}
which  by construction has to be a trace with respect to $\star$
and is $\ad (\A^\infty)$-invariant.
Moreover, by the global symbol calculus for pseudodifferential operators
\cite{Wid:CSCPO,P1} the following formula is satisfied as well:
\begin{equation}
\label{Eq:OpTr}
  \Tr^\text{\tiny \rm Op} (q) = \frac{1}{\hbar^{\dim X}}
  \int_{T^*G_0} q (-,\hbar) \, \omega^{\dim X} .
\end{equation}
Finally, we obtain a trace $\Tr^\text{\tiny \rm Op}_G$ on the algebra
$(\A^{- \infty}_{G,\text{\tiny \rm cpt}} (G_0),\star)$ of
invariant asymptotic symbols as follows:
\begin{equation}
\label{Eq:InvOpTr}
  \Tr^\text{\tiny \rm Op}_G (q) :=
  \Tr_{L^2} \big( \pi_G ( \Op \iota_\hbar (q(-,\hbar))) \pi_G \big).
\end{equation}